\documentclass{article}

\begin{document}

\title{Extremal metrics on toric surfaces, I}
\author{S. K. Donaldson}
\date{\today}
\maketitle

\newtheorem{thm}{Theorem}
\newtheorem{rmk}{Remark}
\newtheorem{prop}{Proposition}
\newtheorem{defn}{Definition}
\newtheorem{case}{Case}
\newtheorem{condn}{Condition}
\newtheorem{cor}{Corollary}
\newtheorem{lem}{Lemma}
\newcommand{\tP}{\tilde{P}}
\newcommand{\tA}{\tilde{A}}
\newcommand{\tsigma}{\tilde{\sigma}}
\newcommand{\tU}{\tilde{U}}
\newcommand{\oP}{\overline{P}}
\newcommand{\tu}{\tilde{u}}
\newcommand{\tx}{\tilde{x}}
\newcommand{\distg}{{\rm Dist}_{g}}
\newcommand{\distEuc}{{\rm Dist}_{{\rm Euc}}}
\newcommand{\End}{{\rm End}}
\newcommand{\uA}{\underline{A}}
\newcommand{\thalf}{\frac{1}{2}}
\newcommand{\tquart}{\frac{1}{4}}
\newcommand{\tthird}{\frac{1}{3}}
\newcommand{\cE}{{\cal E}}
\newcommand{\bR}{{\bf{R}}}
\newcommand{\bC}{{\bf{C}}}
\newcommand{\bP}{{\bf{P}}}
\newcommand{\bZ}{{\bf Z}}
\newcommand{\utheta}{\underline{\theta}}
\newcommand{\Vol}{{\rm Vol}}
\newcommand{\Av}{{\rm Av}}
\newcommand{\db}{\overline{\partial}}
\newcommand{\ux}{\underline{x}}
\newcommand{\oz}{\overline{z}}
\newcommand{\ur}{\underline{r}}
\newcommand{\uSigma}{\underline{\Sigma}}
\newcommand{\uR}{\underline{R}}
\newcommand{\uv}{\underline{v}}
\newcommand{\un}{\underline{n}}
\newcommand{\uC}{\underline{C}}
\newcommand{\ukappa}{\underline{\kappa}}
\newcommand{\tkappa}{\tilde{\kappa}}
\newcommand{\tF}{\tilde{F}}
\newcommand{\dtheta}{\frac{\partial}{\partial \theta}}
\newcommand{\dH}{\frac{\partial}{\partial H}}
\newcommand{\dt}{\frac{\partial}{\partial t}}
\newcommand{\Euc}{{\rm Euc}}
\newcommand{\Riem}{{\rm Riem}}
\newcommand{\Ric}{{\rm Ric}}
\newcommand{\cS}{{\cal S}}
\newcommand{\partialone}{\partial_{x_{1}}}
\newcommand{\partialtwo}{\partial_{x_{2}}}
\newcommand{\partialxi}{\partial_{\xi}}
\newcommand{\partialeta}{\partial_{\eta}}
\newcommand{\ad}{{\rm ad}}
\newcommand{\diag}{{\rm diag}}
\newcommand{\normal}{{\rm norm}}
\newcommand{\Area}{{\rm Area}}
\newcommand{\dbd}{\overline{\partial}\partial}
\newcommand{\cL}{{\cal L}}
\newcommand{\cD}{{\cal D}}
\section{Introduction}
This is the first in a series of papers which continue the study in \cite{kn:D1},
\cite{kn:D2} of the Kahler geometry of toric varieties. The purpose of
the present paper is to introduce an analytical condition (the \lq\lq M-condition'')
and show that it controls sequences of extremal
metrics on toric surfaces. To set the scene
for our discussion we consider the following data:
\begin{itemize}
\item an open polygon $P\subset \bR^{2}$, with compact closure $\oP$;
\item a map $\sigma$ which assigns to each edge $E$ of $P$ a strictly positive
weight $\sigma(E)$;
\item a smooth function $A$ on $\oP$.
\end{itemize}
The datum $\sigma$ yields a measure $d\sigma$ on the boundary $\partial P$---on
each edge $E$ we take $d\sigma$ to be a constant multiple of the standard
Lebesgue measure with the constant normalised so that the mass of the edge
is  $\sigma(E)$. Equally, the datum $\sigma$ specifies an affine-linear defining
function $\lambda_{E}$ for each edge $E$, {\it i.e.} the edge lies in the
hyperplane $\lambda_{E}^{-1}(0)$. We choose an inward-pointing normal vector
$v$ at a point of $E$ with 
$$  \vert i_{V}d\mu \vert =d\sigma_{E} $$
where $d\mu$ is the fixed standard area form on $\bR^{2}$ and we specify
$\lambda_{E}$ by the condition that $\nabla_{v}\lambda_{E}=1$.

For a continuous function $f$ on $\oP$ we set 
$$  L_{A,\sigma} f = \int_{\partial P} f d\sigma - \int_{P} A f d\mu. $$
We require our data $(P,\sigma, A)$ to satisfy the condition that $L_{A,\sigma}
f$ vanishes for all 
affine-linear functions $f$---in other words, that $ \partial P$ and $P$ have
the same mass and centre of mass with respect to the measures $d\sigma$ and
$A d\mu$ respectively. Notice that given $P$ and $\sigma$ there is a unique affine-linear
function $A_{\sigma}$ such that $(P,\sigma, A_{\sigma})$ satisfies this requirement.

Now let $u$ be a convex function on $\oP$, smooth in the interior. We say
that $u$ satisfies the {\it Guillemin boundary conditions} if
\begin{itemize}
\item  any  point $x_{0}$ in the interior of an edge $E$  is  contained in
a neighbourhood
$N_{x_{0}}$ on which 
$$  u= \lambda_{E} \log \lambda_{E} + f $$
where $f$ is smooth in $N_{x_{0}}\cap \oP$ and with strictly positive second
derivative on $N_{x_{0}}\cap E$;
\item if $x_{0}$ is a vertex of $P$, the intersection of two edges $E,E'$,
then there is a neighbourhood $N_{x_{0}}$ on which
$$  u=\lambda_{E} \log \lambda_{E} + \lambda_{E'} \log \lambda_{E'} + f$$
where $f$ is smooth in $N_{x_{0}} \cap \oP$. \end{itemize}

(Note that it is these boundary conditions depend on the weights via the affine-linear defining functions. Thus we can extend the concept to unbounded polygons
with specified defining functions.)

With this material in place, we can recall that the basic question we wish  to
address is the existence of a smooth solution $u$ to the fourth order partial differential
equation {\it (Abreu's equation)}
$$  u^{ij}_{ij}= -A, $$
in $P$, satisfying the Guillemin boundary conditions. (Here we use the summation convention, and $u^{ij}$ is the inverse of the Hessian
of $u$. Our general practice is to use upper indices $(x^{1}, x^{2})$ for
the co-ordinates on $\bR^{2}$, although  we switch to lower indices when
this is more convenient.) If such a function $u$ exists, it is an absolute minimum of the functional
 $$  {\cal F}(f)= - \int_{P} \log \det (f_{ij}) + L_{A,\sigma} f, $$
 over all convex functions $f$ on $\oP$, smooth in the interior. In \cite{kn:D1}
 we were lead to conjecture that a solution exists if and only if the linear
 functional has the property that
 $L_{A,\sigma} f\geq 0$ for all convex $f$ having $L^{1}$ boundary values,
  with strict inequality if $f$
 is not affine-linear. We showed in \cite{kn:D1} that this is a necessary
 condition for the existence of a solution and the problem is to establish
 the sufficiency. We will write ${\cal C}(P)$ for the set of pairs $(A,\sigma)$
 which satisfy this positivity condition.

 The motivation for this problem stems from the case when $P$ is a \lq\lq
 Delzant polygon'',
 corresponding to a compact symplectic $4$-manifold $X$ with a torus action.
 Such a polygon comes with a preferred choice of $\sigma$--we will refer
 to the pair $(P,\sigma)$ as a \lq\lq Delzant weighted polygon''.  The
  convex functions $u$ satisfying the Guillemin boundary
 conditions correspond to invariant Kahler metrics on $X$. In general for
 a strictly convex smooth function $u$ on a polygon $P$ we let $g$ be the Riemannian metric
 on $P$ defined by the Hessian $u_{ij}$ and $\hat{g}$ be its extension to
 $P\times \bR^{2}$ given by
 \begin{equation}
    \hat{g}= u_{ij} dx^{i} dx^{j} + u^{ij} d\theta_{i} d\theta_{j}. \end{equation}
    This is a Kahler metric, with Kahler form $dx^{i}d\theta_{i}$, invariant
    under translations in the $\bR^{2}$ variables. In particular $\hat{g}$
    descends to a metric (which we denote by the same symbol) on $P\times
    \bR^{2}/ 2\pi \bZ^{2}$. If the polygon is Delzant then, with the preferred
    choice of $\sigma$, this metric extends to a smooth metric on a compact
    $4$-manifold $X$. 
    The expression $-u^{ij}_{ij}$
 gives one half the  scalar curvature of the metric $\hat{g}$, \cite{kn:Abr}. When $A=A_{\sigma}$ our problem
 is equivalent to the existence of an  {\it extremal Kahler metric} (in the
 given cohomology class) on $X$. In  particular, if it happens that $A_{\sigma}$ is constant ({\it i.e.} if the
centre of mass of $(\partial P, d\sigma)$ coincides with the centre of mass
of $P$) our problem is equivalent to the existence of a constant scalar curvature
Kahler metric. The positivity condition described above is related to algebro-geometric
notions of \lq\lq stability''.

 In \cite{kn:D1} we obtained a rather weak existence result by the variational
 method applied to the functional ${\cal F}$. In the present paper we change
 our approach to the continuity method. In Section 2 we set up the framework
 for this. We show that solutions persist under small perturbations 
 of the data $(P,A,\sigma)$. Given any polygon $P_{1}$ and $(A_{1}, \sigma_{1})\in
 {\cal C}(P_{1})$ we show that there is a path $(P_{t}, A_{t}, \sigma_{t})$ for $t\in
 [0,1]$ such that $(A_{t},\sigma_{t})\in {\cal C}(P_{t})$ for each $t$ and
 a solution to our problem exists when $t=0$. This is rather trivial if one
 allows arbitrary functions $A_{t}$ but we show that if $A_{1}$ is linear
 (respectively, constant) we can arrange that the $A_{t}$ are also linear
 (respectively, constant). Thus in the standard fashion our problem comes
 down to establishing closedness with respect to $t$, that is to say to establishing
 {\it
 a priori} estimates for a solution $u$ in terms of given data $(P,A,\sigma)$.
   
    In \cite{kn:D2} we studied this problem in the interior
 of the polygon and showed that, roughly speaking, singularities cannot develop
 there. The goal of this paper, and its sequels, is to extend these estimates,
 in appropriate form, up to the boundary. Now we will introduce the central
 notion of this paper.  Let
 $u$ be a smooth convex function defined on sone convex set $\Omega\subset
 {\bR}^{n}$ and let $p,q$ be distinct points in $\Omega$. Let $\nu$ be the
 unit vector pointing in the direction from $p$ to $q$. We write
 \begin{equation} V(p,q)= \left(\nabla_{\nu} u\right)(q)-\left(\nabla_{\nu} u\right) (p), \end{equation}
 where $\nabla_{\nu}$ denotes the derivative in the direction $\nu$. Thus
 $V(p,q)$ is positive by the convexity condition.
 Let $I(p,q)$ be the line segment $$I(p,q)= \{\frac{p+q}{2} + t (p-q): -3/2\leq t\leq
 3/2\}. $$ 
 \begin{defn}
 For $M>0$ we say that $u$ satisfies the $M$-condition if for any $p,q$ such
 that $I(p,q)\subset \Omega$ we have $V(p,q)\leq M$.
\end{defn} 
 
It is easy to see that if the domain is a polygon $P$ as above, and if $u$
satisfies Guillemin boundary conditions, then $u$ satisfies the $M$-condition
for some $M$. Our main result is
 \begin{thm}
 Let $(P^{(\alpha)}, \sigma^{(\alpha)}, A^{(\alpha)})$ be a sequence of  data sets converging to $(P,A,\sigma)$. Suppose  that for each $\alpha$ there is a solution $u^{(\alpha)}$ to the problem defined by
 $(P^{(\alpha)}, \sigma^{(\alpha)}, A^{(\alpha)})$. If there is an $M>0$ such that each $u^{(\alpha)}$
 satisfies the $M$-condition then there is a  solution of the problem defined by
  $(P,A,\sigma)$.
 \end{thm}

 While it is crucial for our continuity method that we do not restrict attention
 to Delzant polygons, it is easier to outline the proof of Theorem 1 in this
 special situation. In Section 3 we develop a variety of arguments which ultimately
 show that the $M$- condition gives a lower bound on the injectivity radius
  of the metric on the  4-dimensional manifold, in terms of the maximal size of the
  curvature (see Proposition 10 below).  If the curvature were to become large, in the sequence, then after rescaling we are able to obtain \lq\lq blow up limits'' which have zero scalar curvature.  In the special
situation when we are actually working with compact $4$-manifolds these limits
could be obtained as a consequence of general results in Riemannian geometry
but we give proofs (in Section 4) adapted to our particular circumstances,
in order to handle general polygons and also in order to make the paper self-contained.
Then we show that these blow-up limits do not exist. There are essentially
two cases to consider. In one case we can appeal to a more general theorem
of Anderson, but we also give an independent proof for the particular result we
need. In the other case we use a maximum principle argument, based on a result
which we prove in the Appendix. Thus we conclude, from the nonexistence of
these blow-up limits, that in fact the curvature
was bounded in the sequence, which leads to the desired convergence.

 The upshot of all this is that we can prove the existence conjecture of
 \cite{kn:D1} if we can establish an {\it a priori} M-condition on solutions.
 More precisely, for given data $(P,\sigma, A)$ and a choice of base point
 $p_{0}\in P$ we can define
 $$  \lambda(P,\sigma, A)= \sup \int_{\partial P} f d\sigma, $$
 where the supremum runs over positive convex functions $f$ vanishing at
 $p_{0}$ and with $L_{A,\sigma} f =1$. We showed in \cite{kn:D1} that, for data in ${\cal C}(P)$ this $\lambda(P,\sigma, A)$ is finite and the remaining problem is to show that
 solutions to our problem satisfy an $M$-condition, where $M$ will depend, among other
 things,  on $\lambda(P,\sigma,
 A)$. This will be taken up in the sequels to the present paper (although
 the author envisages that the actual argument will be rather more complicated
 than this outline suggests).

\section{The continuity method}
\subsection{Connectedness}
  For a given polygon  $P$ we have defined ${\cal C}_{P}$ to be the set of
  $(A,\sigma)$ such that $L_{A,\sigma}$ is strictly positive on the non-affine
  convex functions. Clearly ${\cal C}_{P}$ is itself a convex
  set. We now define a  \lq\lq canonical weight function'' $\sigma_{P}$ as
  follows. Let $p_{0}$ be the centre of mass of $P$, with the standard Lebesgue
  measure on ${\bf R}^{2}$ and for each edge $E$ of $P$ let $cE$ be the triangle
  with base $E$ and vertex $p_{0}$. Obviously, up to sets of measure $0$,
  the polygon $P$ is decomposed into a disjoint union of these triangles.
   Now define
  $$    \sigma_{P}(E)= {\rm Area}\ (cE). $$
  To simplify notation, and without  loss of generality,  suppose  $p_{0}=0$.
  Clearly the mass of the boundary, in the measure $d\sigma_{P}$, is the
  same as the area of $P$. Further, if $q,q'$ are the endpoints of an edge
  $E$ the centre of mass of $cE$ is $\frac{1}{3}(q+q')$ while the centre
  of mass of $E$ is $\frac{1}{2} (q+q')$. Summing over the edges it follows   that the centre of mass of $\partial P$ is also at $0$. Hence the linear
  function $A_{\sigma_{P}}$ associated to these canonical weights is the
  constant function $1$.
  \begin{lem} The pair $(\sigma_{P}, 1)$ is in ${\cal C}_{P}$.
  \end{lem}
  This is essentially  a result of Zhou and Zhu, (Thm. 0.1 of \cite{kn:ZZ}), but since the proof is very simple we include it here. Take standard
   polar co-ordinates $(r,\theta)$ on ${\bf R}^{2}$. By elementary calculus
   one finds that the measure $d\sigma_{P}$ is given by the $1$-form $\frac{1}{2} r^{2}
   d\theta$, restricted to the boundary. Let $f$ be a convex function on the closure of $P$. Since 
   $L_{(\sigma_{P},1)} (f)$ is unchanged  by the addition of an affine-linear
   function, we can suppose  without
   loss of generality that $f$ achieves its minimum value at the origin,
   and that the minimum value is zero.  Now let the boundary be given by
   the equation $r=R(\theta)$. Then we have, by  convexity,
   $$   f(r,\theta)\leq \frac{r}{R(\theta)} f(R(\theta), \theta). $$
   Thus $$ \int_{P} f d\mu= \int_{0}^{2\pi}\int_{0}^{R(\theta)} f(r,\theta) r dr d\theta\leq \int_{0}^{2\pi} \int_{0}^{R(\theta)} \frac{r^{2}}{R(\theta)}
f(R)
dr d\theta. $$
   Integrating with respect to $r$;
   $$  \int_{P} f d\mu \leq \frac{1}{3} \int_{0}^{2\pi} f(R(\theta),\theta)
  \  R(\theta)^{2} d\theta , $$
   whereas
   $$ \int_{\partial P} f d\sigma_{P} = \frac{1}{2} \int_{0}^{2\pi} f(R(\theta),\theta)
   \ R(\theta)^{2} d\theta. $$
   So $$ L_{(\sigma_{P},1)} (f) =\frac{1}{6} \int_{0}^{2\pi} f(R(\theta),\theta)
   \ R(\theta)^{2} d\theta, $$
   and this is clearly strictly positive if $f$ is not identically zero.
   The argument extends immediately to the case when $f$ only has $L^{1}$
   boundary values. 
   
   We define the notion of a \lq\lq continuous path of polygons'' $P_{t}$ in the
   obvious way: the polygons should have the same number of edges and the
   vertices should vary continuously. Similarly, there is an obvious definition
   of a continuous $1$-parameter family of data sets $(\sigma_{t}, A_{t})$
   corresponding to $P_{t}$. 
   \begin{prop}
   Let $P_{t}\ , t\in[0,1]$ be a continuous path of polygons and suppose
   we have $(\sigma_{0}, A_{0})\in {\cal C}(P_{0}), (\sigma_{1}, A_{1})
   \in {\cal C}(P_{1})$. Then these can be joined by  a continuous $1$-parameter family with   $(\sigma_{t}, A_{t})\in {\cal C}(P_{t})$. If $A_{0},A_{1}$
are affine-linear we can suppose that each $A_{t}$ is affine-linear, and
if $A_{0}, A_{1}$ are constant we can suppose that each $A_{t}$ is constant.
\end{prop}
   First, if $\sigma_{0}= \sigma_{P_{0}},\sigma_{1}= \sigma_{P_{1}}, A_{0}=1,
   A_{1}=1$ we can take $\sigma_{t}=\sigma_{P_{t}}\ , \ A_{t}=1$ for all $t\in [0,1]$.  These lie in ${\cal C}(P_{t})$ by the preceding lemma,
and obviously form a continuous family.   Now, by composing paths, we can
reduce to the case when $P_{1}=P_{0}$ and $\sigma_{0}= \sigma_{P_{0}}\ ,\
A_{0}=1$. Here we just use the linear interpolation, applying the convexity of
${\cal C}(P_{0})$. If $A_{1}$ is affine-linear (respectively constant)
then each $A_{t}$ will be affine-linear (respectively constant), and the
proof is complete. 

\subsection{Openness}

Let $P_{t}\ , t\in [0,1]$ be a continuous $1$-parameter family of polygons and $\sigma_{t}$
a $1$-parameter family of weights. Each edge $E$ of $P_{0}$ varies in
a $1$-parameter family $E(t)$ of edges and we have affine-linear defining
functions $\lambda_{E(t)}:\bR^{2}\rightarrow \bR$. We can choose a continuous
$1$-parameter family of diffeomorphisms $\chi_{t}:P_{0}\rightarrow P_{t}$
such that, near to  each edge $E$,
$$    \lambda_{E(t)} \circ \chi_{t}= \lambda_{E}. $$
(This implies that $\chi_{t}$ is affine-linear  near each vertex of $P_{0}$.)
Then, for small $t$,  a function $u_{t}$ on $P_{t}$ satisfies the Guillemin boundary conditions for $(P_{t},\sigma_{t})$ if and only if $\tu_{t}= u_{t}\circ \chi_{t}$ satisfies the
boundary conditions for $(P_{0}, \sigma_{0})$. In a $1$-parameter family,
we say that $u_{t}$ varies continuously with $t$ if the functions $\tu_{t}-u_{0}$
(which are smooth functions on $P_{0}$) are continuous in $t$, along with
all their multiple derivatives.

In this subsection we prove
\begin{prop}
Let $(P_{t},\sigma_{t}, A_{t})$ be a continuous $1$-parameter family of data
and suppose a solution $u_{0}$ to our problem exists when $t=0$. Then for small $t$
there is a solution $u_{t}$, and $u_{t}$ varies continuously
with $t$.
\end{prop}

Of course, this will be proved by linearising and  applying  the implicit function theorem. On the face of it, this might seem a substantial task,
in view of the singular behaviour of the solutions required by the boundary
conditions, but we will explain that the superficial technical difficulties
evaporate when the problem is set up in a suitable way.

We begin by reviewing the relation between complex and symplectic co-ordinates
in this theory, and the role of the Legendre transform. In this Subsection
it will be more convenient to use lower indices $x_{1}, x_{2}$ for our co-ordinates
on the plane. Consider a  convex
function $u$ on a convex open subset $U$  of  $[0,\infty)^{2}\subset \bR^{2}$ which
satisfies Guillemin boundary conditions along the  intersection of $U$ with
the axes, so  
$$ u= x_{1} \log x_{1} + x_{2} \log x_{2} -x_{1}-x_{2}+ f(x_{1}, x_{2}), $$
where $f$ is smooth on  $U$. We suppose that the derivative $\nabla u$ maps
the set  $U\cap (0,\infty)^{2}$ {\it onto} the dual space, in which case the convexity condition implies that it is 
a diffeomorphism. Then the Legendre transform $\phi(\xi_{1}, \xi_{2})$
is defined on the dual space by the formulas
$$  \xi_{a}= \log x_{a} + \frac{\partial f}{\partial x_{a}}, $$
and $$\phi(\xi_{1}, \xi_{2})= \sum x_{a} \xi_{a} - u(x_{1}, x_{2})=-f(x_{1},
x_{2})+  \sum
x_{a}\frac{\partial f}{\partial x_{a}}. $$
The basic fact that we need is that there is a 1-1 correspondence between pairs $(u,U)$ as above and smooth
 $S^{1}\times S^{1}$-invariant functions $\Phi$ on $\bC^{2}$
  with
$i\db \partial \Phi >0$. This is  given by
$$  \Phi(z_{1}, z_{2}) = \phi(\log \vert z_{1}\vert^{2}, \log \vert z_{2}\vert^{2})
. $$  
Further, if a family $u_{t}$ varies continously with respect to an additional parameter
(in the sense of $C^{\infty}$ convergence of the functions $f_{t}$ on compact subsets of their domains) then the transforms $\Phi_{t}$ vary continously
in $t$ (in the sense of $C^{\infty}$ convergence on compact subsets of $\bC^{2}$).

Now let $(P,\sigma)$ be a weighted polygon and $q$ be a vertex of $P$; the
intersection of two edges $E,E'$. The linear parts of the functions $\lambda_{E},\lambda_{E'}$
give a preferred set of linear coordinates on $\bR^{2}$. If $q'$ is another
vertex the two sets of coordinates differ by an element $G(q,q')\in GL(2,\bR)$.
We next review the \lq\lq standard'' case when all the $G(q,q')$ lie in $GL(2,\bZ)$,
{\it i.e.} when $(P,\sigma)$ is a \lq\lq Delzant'' weighted polygon. In this
case we construct a complex surface $X^{\bC}$ from the data in the following
way. For each vertex $q$ we take a copy $\bC^{2}_{q}$ of $\bC^{2}$ and we
identify points using the $G(q,q')$ acting multiplicatively on the open subsets
$(\bC^{*}_{q})^{2}\equiv (\bC^{*})^{2}$.  Thus if 
$G(q,q')= \left(\begin{array}{cc} a & b\\ c&d \end{array}\right)$ we identify $(z_{1}, z_{2})\in \bC^{2}_{q}$ with $(z'_{1}, z'_{2})\in
\bC^{2}_{q'}$ where
$$ z'_{1} = z_{1}^{a} z_{2}^{b}\ , \ z'_{2}= z_{1}^{c} z_{2}^{d}. $$
In this way we get a complex surface $X^{\bC}$, with a $(\bC^{*})^{2}$-action,
containing an open dense orbit $X^{\bC}_{0}$ which is identified with each of the $(\bC^{*}_{q})^{2}$. We denote the quotient space $X^{\bC}/(S^{1}
\times S^{1})$ by $X$. 
 Any point $v$
in $\bR^{2}$ defines a map $\chi_{v}: X^{\bC}_{0}\rightarrow \bR^{+}$. In
the chart $(\bC^{*}_{q})^{2}$ this is given by $(z_{1}, z_{2})\mapsto \vert
z_{1}\vert^{\alpha} \vert z_{2}\vert^{\beta}$, where $v$ has components $(\alpha,\beta)$
in the coordinates $\lambda_{E}, \lambda_{E'}$. 
    Suppose we have a function $u$ on $P$ which satisfies Guillemin boundary
    conditions. For each vertex $q$ we translate to make $q$ the origin,
    and identify $\oP$ with a convex subset of $[0,\infty)^{2}$ using the
    maps $\lambda_{E}, \lambda_{E'}$.   We take the Legendre transform $\phi_{q} $  and pass to logarithmic coordinates to obtain
   a smooth function $\Phi_{q}$ on $\bC^{2}_{q}$. This yields  a collection
   of functions $(\Phi_{q})$ in our charts which satisfy:
   \begin{enumerate}
   \item $i\dbd \Phi_{q}>0$,
   \item $\Phi_{q}$ is invariant under the action of $S^{1}\times S^{1}$,
   \item $\Phi_{q}-\Psi_{q'}= \log \chi_{q-q'}$ on $X^{\bC}_{0}$.
   \end{enumerate}
    Conversely, given such a collection $\Phi_{q}$, we can recover $u$, up
    to the addition of an affine linear function on $P$. Further, the derivative
    of $\phi_{q}$ defines a homeomorphism from $X= X^{\bC}/S^{1}\times S^{1}$
    to $\oP$. 
    
    Next we move on to the case of a general weighted polygon $(P,\sigma)$.
    While we cannot construct a space $X^{\bC}$, we will see that most of the ideas
    above extend. We define a space $X$ by taking for each vertex $q$ a copy
   $[0,\infty)_{q}^{2}$ of $[0,\infty)^{2}$ and identify $(r_{1}, r_{2})$ in $(0,\infty)_{q}^{2}$
   with $(r_{1}^{a} r_{2}^{b}, r_{1}^{c} r_{2}^{d})$ in $(0,\infty)_{q'}^{2}$.
 Of course we can identify $[0,\infty)_{q}^{2}$ with a quotient of  $\bC^{2}_{q}$
 by $S^{1} \times S^{1}$.   This space $X$ has a dense open subset $X_{0}$ on which there are maps
   $\chi_{v}:X_{0}\rightarrow (0,\infty)$. A function $u$ on $P$ satisfying
   Guillemin boundary conditions again yields a collection of functions $\Psi_{q}$
   on $\bC^{2}_{q}$, with the same properties (1), (2), (3) as before, and
   $u$ defines a homeomorphism from $X$ to $\oP$.
   
    Here we digress to consider a general situation. Suppose we have a compact
    topological space $Z$ which is covered by open \lq\lq charts'' $Z_{\alpha}\subset
    Z$. Suppose that for each $\alpha$ there is a homeomorphism from $Z_{\alpha}$
    to $B_{\alpha}/G_{\alpha}$, where $B_{\alpha}$ is the unit ball in some
    Euclidean space and $G_{\alpha}$ is a compact Lie group, acting isometrically
    on the Euclidean space. We suppose we have sheaves $\cL^{p}_{k}$ on $Z$     which
    restrict, in the charts, to the $G_{\alpha}$-invariant locally $L^{p}_{k}$
    functions on the Euclidean spaces (including  $k=\infty$, with  the
    obvious interpretation).   In the case when the $G_{\alpha}$
    are finite groups this is essentially the notion of an orbifold, but
    as far as the author knows there is not a standard terminology for the
    general situation.    The usual machinery of global analysis transfers
    without difficulty to this situation. Thus if we suppose we have a local
    linear operator $D$  taking functions (say) on $Z$ to functions on $Z$, given in the
    charts by a collection of $G_{\alpha}$-equivariant elliptic differential
    operators we can reproduce all the results of the Fredholm alternative,
    invertibility on Sobolev spaces etc. Similarly, for nonlinear operators
    we can apply the usual implicit function theorem arguments, and we will
    not take the space to formalise this further. 
    
    The point of the preceding remarks is that the space $X$ is equipped
    with exactly this kind of structure. It is covered by open sets which
    are identified with quotients $\bC^{2}_{q}/(S^{1} \times S^{1})$, and
    it is easy to see that there are unique sheaves $\cL^{p}_{k}$ as
    above. Thus, while a general weighted polygon does not define a complex
    surface $X^{\bC}$, it does define a space $X$ in which we can apply the
    standard analytical machinery. If we fix a function $u$, and hence an
    identification between $X$ and $\oP$, one easily shows that the \lq\lq smooth'' functions $\cL^{p}_{\infty}$
    on $X$ are identified with the smooth functions on the manifold with
    corners $\oP$, but the situation for general $p,k$ is not so clear and
    in any case we can avoid this issue by working systematically in the
    equivariant charts.

    With all these preliminaries in place, we move on to  our deformation
    problem. First consider the case where we fix the data $(P,\sigma)$ and
    vary the function $A$. Of course we need to stay within the class where
    the mass and centre of mass of $(P,Ad\mu)$ agree with those of $(\partial
    P, \sigma)$.  Working in a chart $\bC^{2}_{q}$, we are in the standard
    situation, considering the
    scalar curvature $S(\Phi)$ of the metric determined by a Kahler potential
    $\Phi$, with $\Phi=\Phi_{q}$. It is well known that this is a nonlinear elliptic differential
    operator. The linearisation has the form
    \begin{equation}  S(\Phi + \eta) = S(\Phi) + \cD^{*}\cD(\eta) + \nabla S. \nabla \eta
    + O(\eta^{2}). \end{equation}    Here $\cD$ is the Lichnerowicz operator $\db_{T} \nabla$, where $\db_{T}$
    is the $\db$-operator on vector fields, and $\cD^{*}$ is the formal adjoint.
     However there is a subtlety here, because the equation we want to solve
     is $S(\Phi)= A$ and while $A$ is a prescribed function on the polygon
     $P$ the identification between $X$ and $\oP$ also depends on $\Phi$,
     so schematically we have an equation $S(\Phi)= A(\Phi)$. Simple calculations
     show that the dependence of $A$ on $\Phi$ precisely cancels out the
     \lq\lq extra'' term in (3). In other words, if we vary our function $A$
     on $P$ to $A+\alpha$ then the linearisation of the equation in the chart
     $\bC^{2}_{q}$ is just
     $\cD^{*}\cD \eta= \alpha$, where $\alpha$ is regarded as a function
     on $\bC^{2}_{q}$ via the identification furnished by $\Phi_{q}$. This
     is rather clear from the \lq\lq moment map'' point of view (compare
     the discussion in \cite{kn:D1}), and we will not take more space to discuss the calculations
     here. The upshot is that we can solve the nonlinear equation, for small
     variations of $A$, provided we avoid the obstructions from the cokernel
     of the linearisation $\cD^{*}\cD$, which is the same as the kernel of
     $\cD$. But this kernel consists exactly of the pull-back of the affine-linear
     functions on $\oP$ and the constraint is just that the mass and centre
     of mass of $\alpha$ vanish,  which is true by hypothesis.

    The case where we deform the data $(P,\sigma)$ is a little more complicated.
    Consider a 1-parameter family $(P_{t}, \sigma_{t})$ of small deformations
    of $(P_{0}, \sigma_{0})$ (in reality the nature of the parameter space
    is irrelevant). Choose a family of diffeomorphisms $\chi_{t}$ as above,
    and let $u_{t}= u_{t} \circ \chi_{t}^{-1}$. Then $u_{t}$ is convex on
    $P_{t}$ and satisfies Guillemin boundary conditions, for small $t$. 
  Fix a vertex $q$ of $P_{0}$ where edges $E,E'$
    meet. There is no loss in supposing that $q$ is the origin and that $\lambda_{E,0},
    \lambda_{E',0}$ are the standard coordinate functions $(x_{1}, x_{2})$.
    The chart $\bC^{2}_{q}$ is regarded as a fixed space, independent of
    $t$, and for small $t$ we have a function $\Phi_{q,t}$ on $\bC^{2}_{q}$
    obtained from the Legendre transform of $u_{t}$. Unwinding the definitions
    $\Phi_{q,t}(z_{1}, z_{2})= \phi_{q,t}(\log\vert z_{1} \vert, \log \vert
    z_{2} \vert)$ where $\psi_{q,t}$ is the Legendre transform of a function
    $u^{*}_{t}$ on a convex set $U_{t} \subset [0,\infty^{2}$. The function
    $u^{*}_{t}$ has the form $u\circ W_{t}^{-1}$, where $W_{t}$    is a diffeomorphism from $U_{0}$ to $U(t)$ which we write as $\tilde{W}(x_{1},
   x_{2})= \tilde{x}_{1}(x_{1}, x_{2}), \tilde{x}_{2}(x_{1}, x_{2})$. This
   diffeomorphism has the property that $\tilde{x}_{i}= x_{i}$ when $x_{i}$
   is small, in particular it is the identity in a neighbourhood
 of the origin, so $u^{*}_{t}=u $ near the origin. It is clear
 then that $u^{*}_{t}$ converges to $u $ as $t\rightarrow
 0$, in $C^{\infty}$ on compact sets. Thus the corresponding functions $\Phi_{q,t}$
 converge, in $C^{\infty}$ on compact sets by the remarks above. 
 
 The conclusion of the discussion above is the following. Let $X_{t}$ be
 the space associated with $(P_{t},\sigma_{t})$.  For small $t$ and
 each vertex $q$ of $P_{0}$ we have an atlas of \lq \lq charts''  
 $$\pi _{q,t}: \bC^{2}_{q} \rightarrow X_{t}  $$
  covering $X_{t}$. In these charts the equation $S(\Phi)=A$ we want to solve
  is given by a continuously varying family of nonlinear elliptic PDE for
  invariant functions. Thus, as before, we can adapt the usual theory from
  the manifold case to construct solutions.

   \subsection{A starting point}

It is clear that any two plane polygons with the same number of edges can be
joined by a continuous path. The next issue we need to address is the existence
of {\it some} data set for which a solution to our problem exists. This is
trivial if we allow arbitrary functions $A$, but for later developments we
want to be able to restrict to the cases where $A$ is constant.

\begin{prop}
For each $r\geq 3$ there is a polygon $P$ with $r$ vertices and a set of
weights $\sigma$ such that there is a solution to our problem for the data $(P,\sigma,
1)$. 
\end{prop}

One possible approach to this is to consider the canonical weights $\sigma_{P}$
associated to any polygon $P$. In this case a solution to the constant scalar
curvature equation must actually satisfy a second order equation of Monge- Ampere type, corresponding (in the local complex differential geometry) to
a Kahler-Einstein metric. This equation, expressed on $P$, is
 $$\log \det (u_{ij})= u- x^{i} u_{i}. $$
  Then one can hope to
extend the proof by Wang and Zhu \cite{kn:ZW} of the existence of Kahler-Einstein metrics on toric Fano varieties to the case of a general polygon $P$. However instead we
will outline another approach by adapting arguments of Arezzo and Pacard
\cite{kn:AP1},\cite{kn:AP2}.

Suppose first  that $(P,\sigma)$ is a Delzant weighted polygon, with 
$A_{\sigma}=1$ (which is, in other langauge, the vanishing of the \lq\lq Futaki invariant'').
One way in which this vanishing  condition can occur is if $P$
is symmetrical about the origin under the map $x\mapsto -x$, and for simplicity
let us suppose that this is the case. Suppose we know that the polarised variety $X$ corresponding
to $P$ admits a constant scalar curvature metric. Now Arezzo and Pacard
study the following general problem: if we know that a complex surface $Z$
admits a constant scalar metric, find a constant scalar curvature metric on the blow-up $\hat{Z}$ of
$Z$ at some finite set of points $z_{1}, \dots z_{q}$ in $Z$. In this problem there is a positive
real parameter associated to each point: the integral of the class of the
 Kahler form on the corresponding exceptional divisor. Arezzo and Pacard show that
one can find such a metric, for small values of these parameters, modulo obstructions coming from the kernel ${\cal H}$ of
the operator $\cD$ on $Z$. Thus there is a  smooth map $F:[0,\infty)^{q} \rightarrow
{\cal H}$ with $F(0)=0$ and the zeros of $F$ in $(0,\delta)^{q}$ give constant
scalar curvature metrics. Now in our case we choose a pair of vertices $q,-q$ of $P$. These
correspond to points, $Q$ and $-Q$ say, in $X$, which are fixed points of the torus action.
Then  the blow up $\hat{X}$ is another
toric surface. 

The translation of the blow-up construction to the language of polygons is
well-known. Choose coordinates, as in the previous subsection, so that $q$ is the origin and $\lambda_{E},\lambda_{E'}$ are the standard coordinate
functions $x^{1}, x^{2}$. Then for small $\epsilon$ we form a new polygon
by removing the triangle $$\{ (x^{1}, x^{2}): x^{1}> 0, x^{2}> 0 ,
x^{1}+x^{2} \geq \epsilon \}$$
from $P$. This operation corresponds to blowing up the point $Q$, and $\epsilon$
to the blow-up parameter mentioned above. The boundary measure on the new
polygon is fixed as follows. On the portion of the boundary which coincides
with the boundary of $P$ the measure is the same as the original one. On
the \lq\lq new'' piece of boundary, corresponding to $x^{1}+x^{2}=\epsilon$
in the coordinates above, the measure is chosen so that the mass of the new
edge is the same as each of the portions of the original edges which were
removed.
  Of course when we blow up both
points $Q,-Q$ we \lq\lq cut off'' two triangles, one with a vertex at $q$
and one with a vertex at $-q$. If we choose the blow-up parameters to be
equal then the new polygon $P_{\epsilon}$ has the same symmetry under $x\mapsto -x$. In
this situation the obstructions arising from the kernel of $\cD$---i.e. from
the affine-linear functions on $P$, are forced to vanish by the symmetry
and it follows directly from the results of Arrezzo and Pacard that there
is a solution of our problem on $P_{\epsilon}$, for small enough $\epsilon$
and  suitable weights. 

This argument comes close to solving our problem. We can start with the
square, corresponding to the manifold $S^{2}\times S^{2}$ with a standard
constant scalar curvature metric. Then cut off two opposite corners to get
a solution for a hexagon, symmetric about the origin. Then cut off two opposite
corners of this to get a solution for an octagon, and so on. Thus we find
r-gons  admitting  solutions for any {\it even} value of $r$. 

Perhaps this argument can be extended by some elementary trick to cover
odd values of $r$ but, lacking this, we go back to appeal to the core idea underlying
Arezzo and Pacard's construction, adapted to the toric situation. They take
the standard zero scalar curvature \lq\lq Burns metric'' on the blow up of
$\bC^{2}$ at the origin, which is asymptotically Euclidean, scale this by
a small factor and glue  it to the original metric on $Z$ to obtain
an \lq\lq approximate solution'' on the blow up. Then the heart of the matter
is to study the problem of deforming this to a genuine solution, via an implicit
function theorem and analysis of the linearised equation. Just as in the
previous subsection, in the toric case the the space $X^{\bC}$ itself plays
no real role here and everything can be formulated in terms of corresponding
operations on the space $X$, using identical local formulae in our equivariant
charts. Further, also as in the previous subsection, the obstructions to
finding a solution can be completely understood in terms of the centre of
mass of the measure $\sigma$. 

Let $P$ be a polygon with at least $4$ vertices and centre of mass at the
origin. Let $q$ be a vertex of $P$
and let $F,F'$ be two edges of $P$ which do {\it not} contain $q$. Let $\sigma$
be a weight function on $\partial P$ such that $A_{\sigma}=1$ and suppose
that there is a solution $u$ of our problem for the data $P,\sigma$, i.e.
a constant scalar curvature metric. Take two positive real parameters $\lambda, \mu$
and consider the family of weight functions $\sigma(\lambda, \mu)$ on $\partial
P$ with
$$ \sigma_{\lambda,\mu}(F)=\lambda \sigma(F)\ ,\  \sigma_{\lambda,\mu}(F')= \mu
\sigma(F')$$
and with $\sigma_{\lambda,\mu}$ equal to $\sigma$ on all the other edges.
Then the centre of mass of $(\partial P,\sigma_{\lambda,\mu})$ yields a map
from $\bR^{+}\times \bR^{+}$ to $\bR^{2}$, and it is easy to see that the
derivative has  rank $2$ at the point $\lambda=\mu=1$. Now take another small
parameter $\epsilon$ and define a polygon $P_{\epsilon}$ by cutting off 
a small triangle at $q$, using this parameter, in the manner discussed above.
For each $\lambda, \mu$ we get a weight function $\hat{\sigma}_{\lambda,\mu}$
for $P_{\epsilon}$. Let $v(\lambda,\mu,\epsilon)\in \bR^{2}$ be the difference
of the centre of mass of $P_{\epsilon}$ and $(\partial P_{\epsilon}, \hat{\sigma}_{\lambda,\mu})$.
The implicit function theorem implies that there are smooth functions $\lambda(\epsilon),
\mu(\epsilon)$ such that $\lambda(0)=\mu(0)=1$ and 
   $$ v(\lambda(\epsilon), \mu(\epsilon), \epsilon)=0. $$
   (Of course what is involved here is just elementary geometry, and one could
   write these functions down explicitly if desired.) 
   This means that, when $\lambda=\lambda(\epsilon), \mu=\mu(\epsilon)$ the
   data $(P_{\epsilon}, \hat{\sigma}_{\lambda,\mu})$ satisfies the obvious
   necessary condition to have a constant scalar curvature metric, i.e.
   $A_{\hat{\sigma}_{\lambda, \mu}}$ is constant. Adapting the proof of Arrezzo
   and Pacard one can show that there is indeed a solution, for small enough
   $\epsilon$. Using this repeatedly we get $r$-gons admitting solutions
   for all $r\geq 4$. When $r=3$ we can use the standard solution coming
   from the Fubini Study metric on $\bC\bP^{2}$ and thus complete the proof
   of Proposition 3.

\section{Geometric estimates}
\subsection{Riemannian geometry in the polygon}

Throughout this section we consider a function $u$ on a polygon $P\subset
\bR^{2}$ as before, satisfying Guillemin boundary conditions determined by
a weight function $\sigma$. We consider
the Riemannian metric $g$ on $P$ defined by the Hessian $u_{ij}$, along with
its extension $\hat{g}$ to $P\times \bR^{2}$. Then $P$ can be regarded as
a totally geodesic submanifold of $P\times\bR^{2}$. Suppose, momentarily,
that the data $(P,\sigma)$ is Delzant, so corresponds to a genuine $4$-manifold $X^{\bC}$,
a compactification of $P\times T^{2}$. Then there is an isometric involution
of $X^{\bC}$ (given by $\theta_{i}\mapsto -\theta_{i}$) with fixed set a smooth
surface $\Sigma$ which can be obtained by gluing $4$ copies of $\oP$ along
suitable edges, and the metric $g$ extends smoothly to $\Sigma$. It is easy
to see from this that, in any case, the metric $g$ extends to a Riemannian
metric on $\oP$, equipped with a suitable smooth structure (as a 2-manifold
with corners), and that the edges are  geodesics. Thus $\oP$ is geodesically
convex, in that any two points can be joined by a minimal geodesic, and any
geodesic can be extended until it reaches the boundary. A main theme of this
subsection is to relate the Riemannian geometry and the Euclidean geometry
in $P$. We write $\distg$ for the distance function defined by $g$ and $\distEuc$
for the Euclidean distance.
  Recall from \cite{kn:D2}, Sec. 5.2 that the tensor 
  $$F^{ij}_{kl}= u^{ij}_{kl}= \frac{\partial^{2} u^{ij}}{\partial x^{k} \partial
  x^{l}}$$ defined by the
  function $u$ is equivalent to the Riemann curvature tensor of the metric $\hat{g}$. We define
$$ \vert F\vert^{2} = F^{ij}_{kl}F^{ab}_{cd} u_{ia} u_{jb} u^{kc} u^{ld}.
$$
Then the absolute value of the sectional curvatures of $\hat{g}$ are bounded
by $\vert F\vert$. In this section we will explore the interaction between the
$M$-condition and a bound on $\vert F\vert$. A crucial fact that we will
use later in the paper is that if $u_{ij}^{ij}=-A$ then
\begin{equation} \int_{P} \vert F \vert^{2} d\mu_{\Euc}- \int_{P} A^{2} d\mu_{Euc}
\end{equation}
is an invariant of the data $(P,\sigma)$, see \cite{kn:D2}, Corollary 5.

\begin{lem}
Suppose $u$ satisfies the $M$-condition. Let $I$ be a line segment in $\oP$
with mid-point $p$ and let $p'$ be an end point of $I$. Then the Riemannian
length of the segment $pp'$ is at most  $$
 \frac{1}{(\sqrt{2}-1)} \sqrt{M} \sqrt{\vert p-p'\vert_{Euc}}
 .$$ 
\end{lem}

We can suppose that $p'$  is the origin
and that $p$ is $(L,0)$, so $\vert p-p'\vert_{Euc}=L$ and the segment of
the $x^{1}$-axis from $0$ to $2L$ lies in $\oP$. We  apply the definition
of
the $M$-condition to the pair of points $p,q$, where $q=(L/2,0)$. This gives
$$\int_{L/2}^{L} u_{11}(t,0) dt \leq M. $$
The Riemannian length of the straight line segment from $q$ to $p$ is
$$  \int_{L/2}^{L} \sqrt{u_{11}}(t,0) dt $$
which is at most 
$$\sqrt{(L/2)}\left( \int_{L/2}^{L} u_{11}(t,0) dt \right)^{1/2}. $$
hence the Riemannian length of this segment is at most $\sqrt{LM/2}$. Replacing
$p$ by $2^{-r} p$ and summing over $r$ we see that the Riemannian length of the segment from $0$ to $p$ is at most
$$  \sqrt{(ML)} \sum_{r=1}^{\infty} (\frac{1}{\sqrt 2})^{r}, $$
from which the result follows.
        
\begin{cor}
Suppose that $u$ satisfies the $M$ condition and that $p$ is a point of $P$.
Then
$$ \distg(p,\partial P) \leq \frac{1}{\sqrt{2}-1} \sqrt{M} \sqrt{ \distEuc(p,\partial
P)}. $$
\end{cor}

To see this we take $p'$ to be the point on $\partial P$ closest to $p$,
in the Euclidean metric. If $p''=2p-p'$ then the segment $p' p''$ lies in
$\oP$ and we can apply the Lemma above.

Next we derive a crucial result which relates the restriction of $u$ to lines
and the curvature tensor $F$.

         \begin{lem}
            At each point of $P$,
             $$\left(\frac{\partial}{\partial x^{1}}\right)^{2} \left(
            u_{11}^{-1}\right) \leq \vert F\vert.$$
  \end{lem}
  
   One way of approaching this
  is to observe that the restriction of the function $u$ to a slice $\{x_{2}=
  {\rm constant}\}$ represents the metric on a symplectic quotient, and then
  to exploit the fact that curvature increases in holomorphic quotient bundles.
  However we will not explain this further and instead give a direct proof.
  Observe that the quantity $$\left(\frac{\partial}{\partial x^{1}}\right)^{2}
\left(u_{11}^{-1}\right)$$
is unchanged by rescaling $x^{1}$. This means that, by rescaling $x^{1}$
and making a different choice of $x^{2}$, we can suppose that at the point
$p_{0}$ in question $u_{ij}$ is the standard Euclidean tensor. Then the square of
the norm of the curvature tensor at this point is
$$ \vert F \vert^{2}= \sum_{i,j,k,l} \left( u^{ij}_{kl}\right)^{2}, $$
and so  $u^{11}_{11} \leq \vert F \vert$. 
Now, at a general point of $P$ we have
$$  u^{11} = \frac{u_{22}}{ u_{11} u_{22}- u_{21}^{2}}, $$
which gives
$$  u^{11}- u_{11}^{-1} = \frac{u_{12}^{2}}{u_{11}(u_{11} u_{22} - u_{12}^{2})}.
$$
Since $u_{12}$ vanishes at the point $p_{0}$  we have
$$ \left(\frac{\partial}{\partial x^{1}} \right)^{2} \left( u^{11}- u_{11}^{-1}\right)
= 2\frac{ (u_{121})^{2}}{u_{11}^{2} u_{22}} = 2 (u_{121})^{2} \geq 0 $$
at $p_{0}$. So
$$ \left(\frac{\partial}{\partial x^{1}}\right)^{2} u_{11}^{-1} \leq u^{11}_{11} \leq
\vert F\vert. $$

\begin{lem}
Let $p$ be a point of $P$ and $\nu=(\nu^{i})$ a unit vector. Suppose the
segment $\{p+t\nu:-3R\leq t\leq 3R\}$ lies in $P$, that $\vert F\vert \leq
1$ in $P$ and that $u$ satisfies the $M$-condition. Then
$$  u_{ij} \nu^{i} \nu^{j} \leq {\rm Max} \left( \frac{2M}{\pi R}, 2\left( \frac{M}{\pi}\right)^{2}\right).$$
\end{lem}

We can suppose that $\nu$ is the unit vector in the $x^{1}$ direction and
that $p$ is the origin. Let $H(t)= u_{11}(t,0)$. We apply the definition of the $M$-condition to obtain
$$  \int_{-R}^{R} H(t) dt \leq M. $$
By the previous Lemma,
$$   \frac{d^{2}}{dt^{2}} H(t)^{-1}  \leq 1. $$
Suppose $H(0)^{-1}= \epsilon$. Then
$$   H(t)^{-1} \leq \epsilon + C t + \frac{t^{2}}{2}, $$where $C=H'(0)$.
Thus $$ H(t)+H(-t) \geq  \frac{1}{\epsilon+ Ct+ t^{2}/2}+ \frac{1}{\epsilon-Ct
+t^{2}/2} \geq \frac{2}{\epsilon+ t^{2}/2}. $$
This gives
$$ \int_{-R}^{R} H(t) \geq \int_{-R}^{R} \frac{1}{\epsilon +t^{2}/2}= 2 \epsilon^{-1/2}
\int_{0}^{R\epsilon^{-1/2}} \frac{dt}{1+t^{2}/2}. $$
So we have
$$    M\geq \frac{2\sqrt{2}}{\sqrt{\epsilon}}\tan^{-1} \left(\frac{R}{\sqrt{2\epsilon}}\right).$$
Now use the fact that $$ \frac{4}{\pi}\tan^{-1}(z)\geq {\rm Min}(1,z)$$
and a little manipulation to obtain the stated bounds on $\epsilon^{-1}= u_{11}(0,0)$.

\

The  results in the rest of this subsection   depend upon a special feature of the Riemannian metric
$g$, and its relation to the metric $\hat{g}$. Consider the $1$-forms $\epsilon_{i}=dx^{i}$
on $P\times \bR^{2}$. Under the isomorphism between cotangent vectors and
tangent vectors defined by the symplectic form these corresponds to the Killing
fields $\frac{\partial}{\partial \theta_{i}}$. These two Killing fields span
a covariant constant subspace of the tangent space, on the other hand they
are Jacobi fields along any geodesic in $P$. Thus we conclude that the $1$-forms
 $\epsilon_{i}=dx^{i}$
satisfy a Jacobi equation of the schematic  form
$$   \nabla^{2}_{t} \epsilon_{i}= F * \epsilon_{i}, $$
along any geodesic. Expressed in different notation, if $e_{1}, e_{2}$ is
a parallel frame of cotangent vectors along a geodesic and if we write
$\epsilon^{i}=\sum G_{ij} e_{j}$, then the matrix $G(t)$ satisfies an equation of
the form
$$   \frac{d^{2}}{dt^{2}} G = - R G, $$
where $R$ is a symmetric matrix with $\vert \left(R_{ij}\right)\vert \leq \vert F\vert$.  If we
express things in terms of the vector fields $\frac{\partial}{\partial \theta_{i}}$
on the $4$-manifold this almost the same as the standard discussion, as in
\cite{kn:Gr},  of the
Fermi fields associated to the orbits of the isometric action. 

We need a simple  comparison result for Jacobi fields.
\begin{lem}
Suppose that $R(t)$ is a symmetric $k\times k$ matrix-valued function on an interval $(0,a)$
 with $\vert (R(t)) \vert \geq -1$. Suppose
that $\epsilon_{1}(t), \dots \epsilon_{k}(t)$ are $k$-vector solutions of the Jacobi equation
$\epsilon''= - R \epsilon$ which are linearly independent at each point in
the interval and with $(\epsilon'_{i}, \epsilon_{j})=(\epsilon_{i},\epsilon_{j}')$.
Then $ \frac{ \vert \epsilon_{1}(t)\vert}{\sinh t}$ is a decreasing function
of $t$. 
\end{lem} 

The author does not find precisely this result stated in standard textbooks,
so we give a proof, although this follows familiar lines. Fix a point
$t_{0}\in (0,a)$ and consider the derivative of $\vert \epsilon_{1}(t)\vert/\sinh
t$ at $t=t_{0}$. Clearly we can suppose that $\epsilon_{i}(t_{0})$ is the standard
orthonormal frame for the $k$-vectors. In particular $\vert \epsilon_{1}(t_{0})\vert
=1$ and we want to show that
$ (\epsilon'_{1}, \epsilon_{1}) \leq \cosh t_{0}/\sinh t_{0}$ at $t=t_{0}$.
Express $\epsilon_{i}(t)$ in terms of the fixed orthonormal frame by   $\epsilon_{i}= \sum G_{i}^{j} e_{j}$ so $G$ is a solution of the matrix  equation
$G''=- R G$
with $G(t_{0})=1$. Set $S= G' G^{-1}$ so that $S$ satisfies the Ricatti
equation $ S'+ S^{2}= R$. The hypothesis that $(\epsilon'_{i}, \epsilon_{j})=
(\epsilon_{i}, \epsilon'_{j})$ implies that $S(t)$ is symmetric for all $t$.
At $t=t_{0}$ we have $(\epsilon'_{1}, \epsilon_{1})= S_{11}$, the $(1,1)$ entry of the matrix $S$, so it suffices to prove that $S(t_{0})\leq \cosh
t_{0}/\sinh t_{0}$, or equivalently that all the eigenvalues of $S(t_{0})$
are bounded above by $\cosh t_{0}/\sinh t_{0}$. Now each eigenvalue $\lambda(t)$
of $S(t)$ satisfies a scalar Ricatti differential inequality
$$ \lambda'+ \lambda^{2}\leq 1 $$ 
 (see \cite{kn:Gr},\cite{kn:P}: by standard arguments we may ignore the complications
that might occur from multiple eigenvalues). Suppose that $\lambda(t_{0})
>\cosh t_{0}/\sinh t_{0}$ Then we can find $\tau\in (0,t_{0})$ such that
$\lambda(t_{0})= \cosh(t_{0}-\tau)/\sinh(t_{0}-\tau)$. Now the function
$\mu(t)=\cosh(t-\tau)/\sinh(t-\tau)$ satisfies the equation
$ \mu' + \mu^{2} = 1$. So $\lambda'-\lambda^{2}\leq \mu'-\mu^{2}$ in the
interval $(\tau, t_{0}]$ and $\lambda(t_{0})=\mu(t_{0})$. It follows that
$\lambda(t)\geq \mu(t)$ for $t\in (\tau, t_{0})$ and since $\mu(t)\rightarrow
\infty$ as $t$ tends to $\tau$ from above we obtain a contradiction.

Notice that Lemma 5 contains as a special case the familiar Rauch comparison result:
if $\vert\epsilon_{1}\vert\sim t$ as $t\rightarrow 0$ then $\vert \epsilon(t)\vert\leq
 \sinh t$ for 
all $t<a$. Notice also that the hypothesis $(\epsilon'_{i},\epsilon_{j})=(\epsilon_{i},\epsilon'_{j})$
is satisfied in our situation, as one sees by a standard manipulation involving the Lie
brackets of the $\frac{\partial}{\partial \theta_{i}}$.

\begin{lem}
Let $E$ be an edge of the polytope $P$ and suppose that the defining function
$\lambda_{E}$ (determined by $\sigma$) is $x^{1}$. Then if $u$ satisfies
Guillemin boundary conditions and $\vert F\vert\leq
1$ throughout $P$  we have
$$   u^{11}(p) \leq \sinh^{2} \distg(p,E) $$
for any $p$ in $P$.
\end{lem}
To see this we consider a geodesic parametrised by $t\geq 0$, starting at time $0$ on the boundary component $E$. Near the boundary we can describe
the geometry in terms of a $4$-manifold with a group action in the familiar
way. The vector field $\frac{\partial}{\partial \theta_{1}}$ is smooth in
the $4$-manifold and vanishes at $t=0$. The condition that $x^{1}$ is the
normalised
defining function just asserts that this vector field is the generator of
a circle action of period  $2\pi$. It follows that
 $$ \lim_{t\rightarrow 0} t^{-1} \vert \frac{\partial}{\partial \theta_{1}}\vert \leq 1, $$
 (with equality when the geodesic is orthogonal to the edge $E$). . Then,
 by the above,  $\sqrt{u^{11}}= \vert \frac{\partial}{\partial \theta_{1}}\vert \leq \sinh t$ and the result follows. 
\begin{cor}
Let $E$ be an edge of $P$ with defining function $\lambda_{E}$. Then if $\vert
F \vert\leq 1$ we have
$$   \lambda_{E} (p) \leq \cosh (\distg(p,E)) -1 . $$
\end{cor}

Notice that this is an affine-invariant statement. There is no loss in supposing
that, as above, $\lambda_{E}=x^{1}$. Then for a geodesic starting from a point
of $E$, parametrised by arc length, we have
$$   \vert \frac{dx ^{1}}{dt}  \vert\leq \vert dx^{1} \vert_{g} = \sqrt{u^{11}}
\leq \sinh t $$
hence $x^{1} \leq \cosh t -1$.

\begin{lem}
Suppose that $\vert F\vert \leq 1$ and that $p$ is a point in $P$ with $dist_{g}(p,\partial
P)\geq \alpha>0$. Then if $q$ is a point with $\distg(p,q)= d$ we have

 $$\left( u^{ij}(q)\right)\leq \frac{\sinh^{2}(\alpha+d)}{\sinh^{2} \alpha}
\left(u^{ij}(p)\right). $$
If $d<\alpha$ we have
$$ \left(u^{ij}(q)\right)\geq \frac{\sinh^{2}(\alpha - d)}{\sinh^{2} \alpha} \left(u^{ij}(p)\right)
. $$
\end{lem}
(Here the notation $\left( A^{ij}\right) \leq \lambda \left( B^{ij} \right)$
means that for any vector $\nu_{i}$ we have $\nu_{i} \nu_{j} A^{ij} \leq
\lambda \nu_{i} \nu_{j} B^{ij}$.)
To prove the Lemma, observe that it suffices by affine invariance to prove the corresponding inequalities
for the matrix entry $u^{11}=\vert \epsilon_{1}\vert^{2}$. For the first
inequality we consider a minimal geodesic $\gamma$ from $p=\gamma(0)$ to $q=\gamma(d)$ and extend it \lq\lq backwards'' to $t>-\alpha$. Then replacing
$t$ by $t+\alpha$ we are in the situation considered in Lemma 5 and we obtain
$$  \frac{\vert \epsilon_{1}(p)\vert}{\sinh \alpha} \geq \frac{\vert \epsilon_{1}(q)\vert}{\sinh(\alpha
+d)}. $$
For the second inequality we extend the geodesic \lq\lq forwards'' to the interval $[0,\alpha]$
and argue similarly. 

Suppose that $p=(p^{1}, p^{2})$ is a point of $P$ and $r>0$. Put
$$ E_{p,r}= \{ (x^{1}, x^{2})\in \bR^{2}: u_{ij}(p)(x^{i}-p^{i})(x^{j}-p^{j})
\leq r^{2}\}. $$
So $E(p,r)$ is the interior of the ellipse defined by the
parameter $r$ and the quadratic form $u_{ij}(p)$. The Euclidean area of $E(p,r)$
is $\pi r^{2} \det (u_{ij}(p))^{-1/2}$. 

\begin{lem}
Suppose that $\vert F\vert \leq 1$ and that $p$ is a point in $P$ with
$dist_{g}(p,\partial P) \geq \alpha>0$. Then for any $\beta<\alpha$ the
$\beta$-ball in $P$, with respect to the metric $g$ satisfies
$$  E(p,c\beta)\subset   B_{g}(p,\beta) \subset E(p,C\beta), $$
 where $c=\sinh(\alpha-\beta)/\sinh \alpha$ and $C=\sinh(\alpha+\beta)/\sinh
 \alpha$.
In particular, the Euclidean area of the $\beta$ ball for the metric $g$
is bounded below by
$$   {\rm Area}_{{\rm Euc}} B_{g}(p, \beta) \geq \pi c^{2} \beta^{2} \det(u_{ij})(p)^{-1/2}.
$$

\end{lem}

There is no loss in supposing that the matrix $u^{ij}(p)$ is the identity
matrix, so we have to show that the ball $B_{g}(p,\beta)$ defined by the metric $g$ contains
a Euclidean disc of  radius $c\beta$,a nd is contained in a Euclidean disc
of radius $C\beta$.   
We know by Lemma 7  that on the ball $B_{g}(p,\beta)$ we have
$$    c^{2} \leq (u^{ij}) \leq C^{2}$$

Thus $C^{-2} \leq (u_{ij}) \leq c^{-2}$, and the Euclidean
length of a path in $B_{g}(p,\rho)$ is at least $c^{-1}$ times the length
calculated in the metric $g$, and at most $C^{-1}$ times that length. The
second statement immediately tells us that $B_{g}(p,\beta)$ lies in $E(p,C\beta)$.
In the other direction, suppose $q$ is a point in the Euclidean disc
of radius $c\beta$ centred on $p$. We claim that $q$ lies in the (closed) $g$
ball $B_{g}(p,\beta)$. For if not there is point $q'$ in the open line segment
$pq$ such that the distance from $q'$ to $p$ is $\beta$ and the line segment
$pq'$ lies in $B_{g}(p,\beta)$. But the Euclidean length of this line segment
is strictly less than $c\beta$ so the length in the metric $g$ is less than
$\beta$, a contradiction.
      
      \
      
  \subsection{The injectivity radius}

  We continue to consider a convex function $u$, satisfying Guillemin boundary
  conditions, on a polygon $P$, as in the previous subsection. The present
  subsection has two purposes. In one direction we discuss coordinates in
  neighbourhoods of boundary points obtained from geodesic coordinates in
  four dimensions. In another direction, we want to relate these ideas to
  the standard notion of the injectivity radius. Since we will want sometimes
  to work with incomplete manifolds we should clarify our definitions. By
  the statement that \lq\lq the injectivity radius at a point $p$ is at least
  $r$'' we mean that the exponential map at $p$ is defined on tangent vectors
  of length $r$, and yields an embedding of the Euclidean $r$-ball. In fact
  the discussion of the injectivity radius need only enter our main proof
  in a rather minor way, but it is useful to explain how the arguments fit
  into the wider world of Riemannian geometry.
  
  First we consider the vertices.   Let $q$ be
  a vertex of $P$, so  we have an  Riemannian $4$-manifold $X^{c}_{q}$, which
  is not complete. The torus action on $X^{c}_{q}$ gives a constraint on
  the exponential map.
  \begin{lem}
  If $\vert F \vert \leq 1$ in $P$ then the injectivity radius of $X^{c}_{q}$
  at $q$ is at least $\pi/2$.
  \end{lem}
  The exponential map is equivariant with respect to the standard
  torus action on the tangent space at $q$. Suppose the exponential map is
  defined for some $r'<r$ and let $\xi$ be a unit vector in the Lie algebra
  of the torus, corresponding to a vector field $v_{\xi}$ on $X^{c}_{q}$.
  Then the length of the vector field $v_{\xi}$ is bounded below on the boundary
  of the $r'$ ball. However, if the exponential map is not defined on the $r$ ball
  then as we let $r'$ approach its maximal possible value there is some choice
  of $\xi$ such that the length of $v_{\xi}$ goes to zero on the boundary
  (since the corresponding points in $\oP$ must be approaching another edge).
  
  Now suppose that the $r$ ball is not embedded by the exponential map. Then
  there is a nontrivial  geodesic starting and ending at $q$, of length less
  than $\pi$. But the vector fields $\frac{\partial}{\partial \theta_{i}}$
  give Jacobi fields along this geodesic, vanishing at the endpoints. By
  a standard comparison theorem these vector fields must vanish identically
  along the geodesic which means that the initial tangent vector of the geodesic
  is fixed by the torus action. Since there are no such fixed tangent vectors
  we have a contradiction.                           
    
 Now we note a general fact of Riemannian geometry.   

    \begin{lem}
    Let $g_{ij}= \delta_{ij} + \eta_{ij}$ be a  Riemannian metric on the
    Euclidean
     ball $B$ of radius $\pi$ in $\bR^{n}$,  with sectional curvature bounded in absolute value
    by $1$. Suppose that $\vert \left( \eta_{ij} \right)\vert \leq \epsilon \left(\delta_{ij}\right)$,
    for some $\epsilon<1$. Then the injectivity radius at the origin is at
    least $\sqrt{1-\epsilon}$.
    
    \end{lem}
  First, the $g$-distance from the origin to the boundary of the ball $B$ is at least
  $\pi\sqrt{1-\epsilon}$, so the exponential map is defined as stated. Since
  the curvature is less than $1$, we only need to check that there are no
  geodesic loops starting and ending at the origin, of length less than
  $2\sqrt{1-\epsilon}$.     Suppose $\gamma$ is a geodesic loop, of length $L$, and for $s<1$ let $\gamma_{s}$ be the loop $\gamma_{s}(t) = s \gamma(t)$. Then the length of
$\gamma_{s}$ is at most $L'=\sqrt{\frac{1+\epsilon}{1-\epsilon}} L$. For small
$s$ the loop $\gamma_{s}$ can be lifted to a loop over the exponential map.
The
 argument on page 100 of \cite{kn:CE} (proof a Theorem of Klingenberg) shows that this is true for all $s$, provided that $L'<\pi$,
 which will be the case if $L<2\sqrt{1-\epsilon}$. But, as in the argument
 cited, $\gamma$ itself
lifts to a ray under the exponential map, giving a contradiction.  
    
  \
  
  Now we consider an interior point $q$ of the polygon. We can think of this
  as a point in the Riemannian $4$-manifold $P\times \bR^{2}$, with the metric
  $\hat{g}$ and we write $I(q,\hat{g})$ for the injectivity radius at that
  point. We can also consider $q$ as a point in the quotient space $P\times \bR^{2}/\bZ^{2}$ and we write $I'(q,\hat{g})$ for the injectivity radius
there. 

  \begin{lem}
   Suppose $\vert F\vert \leq 1$ in $P$.
  \begin{enumerate}\item For any $\alpha>0$ there is an $i(\alpha)>0$ such that if  $\distg(q,\partial
  P)\geq \alpha$ then $I(q,\hat{g})\geq i(\alpha)$.
  \item If $u$ satisfies an $M$ condition then there is an $i(\alpha,M)>0$
  such that if $\distg(q,\partial P)\geq \alpha$ then $I'(q,\hat{g})\geq
  i'(\alpha,M)$.
  \end{enumerate} 
  \end{lem}
  To prove the first item we apply Lemma 8. We can suppose that the Hessian
  $u_{ij}$ at the point $q$ is the standard form $\delta_{ij}$. Then Lemma
  8 tells us that that the metric $\hat{g}$ is close to Euclidean---in the
  given coordinates $x^{i}, \theta_{j}$---over a
  ball of a definite size determined by $\alpha$. Then we can apply Lemma
  10. To prove the second item we just need to check that the quotient  by
  $\bZ^{2}$ does not create any short loops. Since the metric in the fibre
  direction is given by $u^{ij} d\theta_{i} d\theta_{j}$, this is the same
  as showing  that for any non-zero integer
  vector $\nu_{i}$ the quantity $u^{ij} \nu_{i} \nu_{j}$ is not small. But
  we know, by combining Lemmas 2 and 4, that $u_{ij} \leq C$, where $C$ depends
  on $M,\alpha$. This implies that $u^{ij} \nu_{i} \nu_{j} \geq C^{-1}(\nu_{1}^{2}
  +\nu_{2}^{2}) \geq C^{-1}$.
  
  \
  
  To take stock of our progress so far, consider the case when $P$ corresponds
  to a compact $4$-manifold $X^{c}$. Then Lemmas 9 and 11 give lower bounds
  on the injectivity radius at points of $X^{c}$ which correspond to  either vertices
  or to interior points of $P$. Our remaining task is to consider the points
  which lie on the boundary edges. For this we introduce a numerical invariant
  of a weighted polygon $(P,\sigma)$. Let $q$ be a point in the interior
  of an edge $E$ and let $d$ be the Euclidean distance from $q$ to the end
  points of $E$. Set
  $$  \mu(q)= \min_{E'}\frac{ \lambda_{E'}(q)}{d}, $$
  where $E'$ runs over the set of edges {\it not equal} to $E$. Now let $\mu=\mu_{P,\sigma}$
  be the minimum of $\mu(q)$ over all such boundary points $q$. It is easy
  to see that $\mu_{P,\sigma}>0$. 
  
  For each (open) edge $E$ of $P$ we define 
  a Riemannian $4$-manifold $X^{c}_{E}$ as follows. We choose coordinates
  such that the defining function $\lambda_{E}$ is $x^{1}$ and take the quotient
  of $P\times \bR^{2}$ by $\Lambda$, where $\Lambda$ is the copy of ${\bf
  Z}$ embedded as ${\bf Z}\times \{0\}$ in $\bR^{2}$. This gives a manifold
  with an action of $S^{1} \times \bR$. Then, just as in the
  construction of the manifolds $X^{c}_{q}$ associated to vertices $q$, we
  can adjoin a copy of $E\times \bR$, fixed under the circle action, and
  the metric extends smoothly. If $q$ is a point on the interior of $E$ we write $I(q,\hat{g})$ for the injectivity radius about the corresponding point in $X^{c}_{E}$.
If $(P,\sigma)$ is Delzant we can also consider $q$ as a point in the compact
manifold $X^{c}$ and we write $I'(q,\hat{g})$ for the injectivity radius there.
 
\begin{lem} Suppose that $\vert F\vert \leq 1$ in $P$ and that $u$ satisfies
an $M$ condition.
Then for any $\alpha>0$ there is an $i(\alpha,\mu,M)>0$ such
that $I(q,\hat{g})\geq i(\alpha,\mu,M)$ if the distance in the metric $g$ from $q$
to the set of vertices is at least $\alpha$. If $(P,\sigma)$ is Delzant then
there is an $i'(\alpha,\mu,M)>0$ such that $I'(q,\hat{g})\geq i'(\alpha,\mu,M)$
\end{lem}

Not surprisingly, the proof of this Lemma--for an edge point-- is a combination of the arguments
used in the cases of vertices and interior points. The first thing is to
see that the exponential map at $q$ in $X^{c}_{E}$ is defined on a ball of
a definite size (depending on $\alpha, \mu, M$). This is the same as showing
that the distance in the metric $g$ from $q$ to any other edge $E'$ of $P$ is
not small. But we know by Lemma 2 that the Euclidean distance $d$ from $q$
to the end points of $E$ is not small, hence by the definition of $\mu$,
$\lambda_{E'}(q)$ is bounded
below by a quantity depending on $\mu,M,\alpha$.Then Corollary 2 implies that the distance
in the metric $g$ from $q$ to $E'$ is not too small. The remaining task is
to show, as in the proof of Lemma 9, that there are no short geodesic loops in
$X^{c}_{E}$ starting at $q$. Now there is a circle action on $X^{c}_{E}$
which fixes the point $q$ and the argument used in  the proof of Lemma 9 shows that any short
geodesic loop must lie in the fixed set of the action, which is
$E\times \bR$. The Riemannian metric on $E\times \bR$ is defined by the restriction
of $u$ to $E$. The arguments used in the proof of Lemma 11 apply, in an obvious way, to give
a lower bound on the injectivity radius in $E\times \bR$, so we see that
there are no short geodesic loops and the proof of the lower bound on $I(q,\hat{g})$ is complete.

In the case when $(P,\sigma)$ is Delzant a neighbourhood of $q$ in $X^{c}$
is quotient of $X^{c}_{E}$ by an action of  $\bZ$ and we again we need to
show that this does not create any short loops. This just comes down an upper
bound on the second derivative of $u$ along the edge, which is furnished
by Lemma 4 and the M-condition.

    \begin{prop} Suppose that $(P,\sigma)$ is Delzant, that $u$ satisfies an
    $M$ condition and $\vert F \vert \leq 1$ in $P$. Then there is an $r$, depending
    only on $M$ and $\mu_{P,\sigma}$, such that the injectivity radius of the Riemannian
    $4$-manifold $X^{c}$ is at least $r$.
      \end{prop}
    By applying Lemma 10 it suffices to show that for any $\kappa>0$ there is an $r'$ such that for each point $p$ of $X^{c}$ we can find another point $p'$
such that the injectivity radius at $p'$ is at least $r'$ and the distance
from $p$ to $p'$ is at most $\kappa r'$. If $p$ is close to a vertex we take
$p'$ to be the vertex and use Lemma 9. If $p$ is close to an edge but not close
to any vertex we take $p'$ to be a nearby point on the edge, and use Lemma
12.
If $p$ is not close to any edge we take $p'=p$ and use Lemma 11.

\

 We conclude this section with another simple observation, similar to Lemma
 9, which will be useful
 later.
 \begin{lem}
 Suppose that $q$ is a point on an edge $E$ of $P$ and $\gamma$ is a geodesic
 starting at $q$ which is orthogonal to $E$ at $q$. If $p$ is the point a
 distance $d$ from $q$ along the geodesic, where $d<\pi/2$, then $\distg(p,E)=d$.
 \end{lem}
 In the case when $(P,\sigma)$ is Delzant this is essentially a standard
 result. By the same argument as in Lemma 9, a geodesic segment with
 endpoints on $E$ of length less than $\pi$ must lie in $E$. This means that
 the exponential map on the normal bundle of the $2$-sphere corresponding
 to $E$ is an embedding on vectors of length less than $\pi$, from which
 the assertion follows. The reader can easily check that the proof works
 in just the same way for a general $(P,\sigma)$.

       \section{Convergence of sequences}
       \subsection{Elliptic estimates}
       In this subsection we assemble some results of a rather standard nature;
       the general theme being that the derivatives of the scalar curvature of a Kahler metric
        controls those of the full curvature tensor. Similar, but more sophisticated,
        results are contained in \cite{kn:An}, \cite{kn:TV}.
      
       {\it Throughout this subsection} we suppose that
        $(M,g,J)$ is a Kahler surface with scalar curvature $S$ and let $p$ be a point of $M$. We suppose that the exponential map at $p$
       is defined on the unit ball and for $\rho\leq 1$ let $B_{\rho}$ be the
       $\rho$ ball in centred at $p$.

        We begin with a simple result, which will be the essential thing we
       need for our main argument 
       \begin{prop}
       Suppose that $\vert {\rm Riem}\vert \leq 1$ on $B_{1}$.
       Then for any $\alpha\in (0,1)$ and $\rho<1$ there
       is a Holder bound, for points $p'$ with $d(p,p')\leq \rho$, 
       $$\vert\  \vert \Riem(p') \vert -\vert \Riem(p) \vert\ \vert \leq C_{\alpha,\rho}(
       1+
       \Vert \nabla S \Vert_{L^{\infty}})d(p,p')^{\alpha}. $$
       \end{prop}
        By pulling back the metric we can suppose that the exponential map is an embedding on the unit ball. By a covering argument it suffices to prove
the result for some $\rho$ and then by rescaling we can suppose that $\vert
\Riem \vert$ is as small as we please.

Various approaches to the proof  are possible. We will base or argument
on a general perturbation result for linear elliptic equations. Suppose that
$D_{0}$ is a constant-coefficient first order elliptic operator over $\bR^{n}$
(i.e.
with injective symbol) and $E$ is a perturbation term, defined over the unit
ball, of the form
$$  E(f)= \sum \epsilon_{i} \frac{\partial f}{\partial x_{i}} +  T f. $$
(Here we are considering operators on vector-valued functions, so the coefficients
will be matrices in general.) Fix an exponent $p>1$ and suppose that
\begin{itemize} \item 
$\epsilon_{i}$ are sufficiently small;
\item we have $L^{q}$ bounds on $T$
\end{itemize}
where $q$ and the allowable size of the $\epsilon_{i}$ depend on
$D_{0}$ and $p$. Then by considering $D_{0}+E$ as a perturbation of $D_{0}$
we obtain an elliptic estimate of the form
$$  \Vert f \Vert_{L^{p}_{1}(B_{1/2})} \leq C\left( \Vert (D_{0}+E) f \Vert_{L^{p}(B)}
+ \Vert f \Vert_{L^{p}(B)}\right), $$
where $C$ depends on the $L^{q}$ bounds on the coefficients $T$.
The proof is essentially the same as \cite{kn:GT} Theorem 9.11, together with the remark on page 241.

To apply this we work in geodesic coordinates on our Kahler surface. A bound on the curvature gives a $C^{1}$ bound on the the metric coefficients $g_{ij}$in
these co-ordinates. Since the metric is Kahler the almost-complex structure $J$ is
covariant constant hence, when written as a tensor in these coordinates, the coefficients are also bounded in $C^{1}$.  We use the following identities connecting
the curvature tensors, written in a schematic form
$$  \db \Riem=0\ \ ,\ \   \db^{*} \Riem= \pi(\nabla Ric),$$
$$ \db \Ric=0\ \ ,\ \  \db^{*} \Ric = \pi(\nabla S). $$
Here $\pi$ denotes certain natural contractions on tensors of the appropriate
type. Then we can apply the discussion above to the elliptic operator $\db\oplus
\db^{*}$ defined by the Kahler metric. We express this, in geodesic coordinates,
as a perturbation of the constant coefficient model. When the curvature is
small the relevant terms $\epsilon_{i}, T$ are small in $L^{\infty}$. Now
the general elliptic estimate above yields
$$  \Vert \Riem \Vert_{L^{p}_{1}(B_{1/2})} \leq C (\Vert \nabla S \Vert_{L^{p}}
+1), $$
and we get a $C^{\alpha}$ bound on $\vert\Riem\vert$ from the Sobolev embedding
theorem.  

\

Next we extend this to higher derivatives.
\begin{prop}
With notation as above suppose that $\vert \Riem \vert \leq 1$ on $B_{1}$.
Then for any $l\geq 1$ there are constants $C_{l,\rho}$ such
that $$ \vert \nabla^{l} \Riem \vert \leq C_{l,\rho} (\Vert \nabla^{l+1}
S \Vert_{L^{\infty}(B)} + 1), $$
on $B_{\rho}$.
\end{prop}
We only outline a proof, since this is somewhat standard.
We can apply the perturbation argument as above to the
$\db + \db^{*}$- operator mapping from $L^{p}_{k+1}$ to $L^{p}_{k}$ provided
we know that the coefficients $\epsilon_{i}, T$ are controlled in  $L^{p}_{k}$.
(Here $p$ is chosen sufficiently large.) Since $T$ depends on the first derivatives
of the metric tensor $g$ and the complex structure $J$, in coordinates, we
need $g,J\in L^{p}_{k+1}$. To achieve this we work in harmonic coordinates
\cite{kn:JK},
in which the $L^{p}_{k+1}$ norm of the metric tensor is controlled by the
$L^{p}_{k-1}$ norm of the curvature tensor. Since the tensor $J$ is covariant
constant we also get an $L^{p}_{k+1}$ bound on its representative in these
coordinates.. Now
we bootstrap, starting from the $L^{p}_{1}$ bound on the curvature tensor
which was already obtained in the proof of Proposition 5. In harmonic coordinates we
can consider the $\db \oplus \db^{*}$ operator mapping $L^{p}_{3}$ to $L^{p}_{2}$
and obtain $L^{p}_{3}$ bounds on the curvature tensor, in terms of derivatives
of the scalar curvature, and so on.
 
 \
 
Now consider a more specialised situation in which we have a pair of holomorphic vector
fields  $v_{1}, v_{2}$ on an embedded ball $B_{1}$ in the Kahler manifold $X$. Suppose that the Riemannian gradient of the scalar curvature
can be expressed as  $\nabla S = A_{1} v_{1} + A_{2} v_{2}$ where $A_{1},
A_{2}$ are functions on the manifold. Suppose in turn that all derivatives
of $A_{1}, A_{2}$ can be expressed in a similar way:
$$ \nabla A_{i} = \sum A_{ij} v_{j}. $$
$$ \nabla A_{ij} = \sum A_{ijk} v_{k} , $$
and so on. 
\begin{prop}
In this situation, if $\vert \Riem  \vert \leq 1$ on $B_{1}$ then we have
$\vert \nabla^{l} \Riem \vert \leq C_{l,\rho}$
on $B(\rho)$, where $C_{l,\rho}$ depends on the $L^{\infty}$ norms of the
vector fields $v_{1}, v_{2}$ and the functions $A_{i_{1}\dots i_{k}}$ over
the ball $B_{1}$, for $k\leq l+1$.
\end{prop}

To prove this we exploit the first order elliptic equation $\db v_{i}=0$
for the vector fields and build this into our bootstrapping argument.
First, the $L^{\infty}$ norm of $\nabla S$ is obviously controlled by the
$L^{\infty}$ norms of $A_{i}, v_{i}$. So in harmonic coordinates we control
the $L^{p}_{3}$ norm of the metric and obtain elliptic estimates for the
$\db$-operator mapping $L^{p}_{3}$ to $L^{p}_{2}$ and we get an $L^{p}_{3}$
bound on $v_{i}$. Now we can write
$$ \nabla^{2} S = \sum A_{ij} v_{i} \otimes v_{j} + \sum A_{i} \nabla v_{i}$$
and we get an $L^{p}_{2}$ bound on $\nabla^{2} S$ and so on.

\subsection{Bounded curvature}

Now we show that to prove Theorem 1 it suffices to bound the curvature tensors of
the solutions.
\begin{prop}
Suppose that $(P^{(\alpha)},\sigma^{(\alpha)}, A^{(\alpha)})$ are data-sets converging to a limit
$(P,\sigma, A)$ and that $u^{(\alpha)}$ are solutions.
If there are fixed $M,K$ such that $u^{(\alpha)}$ satisfies the $M$-condition and $\vert F(u^{(\alpha)})\vert \leq K$, for all $\alpha$,  then there
is a solution $u^{(\infty)}$ for the data $(P, \sigma, A)$.
\end{prop}
Of course, the solution $u^{(\infty)}$ will be obtained as a limit of the $u^{(\alpha)}$,
provided that these are suitably normalised with respect to the addition
of affine-linear functions. Although
the domains of definition $P^{(\alpha)}$ are different, it obviously makes sense
to talk about a subsequence of the $u^{(\alpha)}$ converging on compact subsets of $P$, and this
is what we show first (In fact we already have this interior convergence from the
results of \cite{kn:D2}---without assuming the curvature bound--- but we will give an independent argument since it will be pave the way for the proofs in 4.4 below.) To simplify
the presentation we just consider the case when the $P^{(\alpha)}, \sigma^{(\alpha)}$ are
all the same $(P,\sigma)$ and only $A^{(\alpha)}$ varies with $\alpha$. The reader will easily see
that the general case is not essentially different. We simplify notation
by sometimes writing $u$ and $A$ for $u^{(\alpha)}$ and $A^{(\alpha)}$.

By Lemma 2, there is some fixed $D$ such that for any  point $p$ in $P$ there is
a vertex $q$ such that the Riemannian distance from $p$ to $q$ is less than
$D$. Then Lemma 6 gives a universal bound
$$   u^{ij} \leq C.$$
On the other hand  Lemma 4 gives a bound
$$  u_{ij} \leq C/ d_{Euc}, $$
where $d_{Euc}$ is the Euclidean distance to the boundary of $P$. So we deduce
that $u_{ij}$ is bounded above and below on compact subsets of the interior.
On such sets the definition of the curvature tensor $u^{ij}_{kl}$ immediately
gives a $C^{2}$ bound on the $u_{ij}$, so we can suppose that the $u_{ij}$
converge in $C^{3,\alpha}$. From this it is entirely straightforward to deduce
the $C^{\infty}$ convergence, on compact subsets of $P$. Thus the essential issue is to show that the limit
satisfies the Guillemin boundary conditions. To see this, fix a point $q$
on the boundary of $P$.   There are
two cases to consider, either $q$ is a vertex or lies on the interior of
an edge $E$. 

{\it Case 1: $q$ is a vertex}

The function $u= u^{(\alpha)}$ defines an $S^{1} \times S^{1}$-invariant metric on $X_{q}\cong \bC^{2}$. By Lemma 9, the geodesic ball of some fixed small radius about the origin is embedded. This geodesic ball maps to neighbourhood of $q$ in $\oP$ which is
contained in a Euclidean neighbourhood of one fixed size, and contains a
Euclidean neighbourhood of another fixed size. We are in the framework of
Proposition 7, with $v_{i}= I \frac{\partial}{\partial \theta_{i}}$ and $A_{i}=\frac{\partial
A}{\partial x^{i}}, A_{ij}= \frac{\partial^{2}}{\partial x^{i} \partial x^{j}}$
and so on. Thus the norm of $v_{i}$ in the Riemannian metric is $\sqrt{u^{11}}$
and this is bounded.  Similarly for $v_{2}$. All the derivatives of $A$ are
bounded so we can apply Proposition 7
to deduce that all covariant derivatives of
the curvature tensor are bounded in this ball. We pass to geodesic coordinates
in which we have data $(g^{(\alpha)}, J^{(\alpha)})$. Then in these geodesic coordinates all
derivatives of the metric tensors are bounded and we can suppose that the
metrics converge in $C^{\infty}$, likewise for the complex structures since
these are covariant constant.  The limit is a
smooth Kahler metric $(g^{\infty},J^{\infty})$ on a
small ball in $\bR^{4}$, invariant under the fixed, standard, action of $S^{1}\times
S^{1}$. For each $\alpha$, the functions $x^{1}_{\alpha}, x^{2}_{\alpha}$ which map the ball to  neighbourhoods
of $q$ in $\oP$ are characterised as moment maps for the action with respect
to the symplectic forms $\omega^{(\alpha)}$ determined by $(g^{(\alpha)}, J^{(\alpha)})$. It follows
that these also converge. By Guillemin's analysis of the structure of invariant
Kahler metrics we know that the limit $(g^{\infty},J^{\infty})$ corresponds to a function $u^{\infty}$
on a neighborhood of $q$ in $\oP$, satisfying Guillemin boundary conditions,
and it is clear from the convergence of the data $(g^{(\alpha)}, J^{(\alpha)}, x^{1}_{\alpha},
x^{2}_{\alpha })$ that the second derivative of this coincides with the limit we have already found on the interior. Thus we see that this interior limit satisfies Guillemin boundary conditions in a neighbourhood
of the vertex $q$.

\

{\it Case 2: $q$ is in the interior of an edge}

We suppose that $P$ is defined near $q$ by the equation $x^{1}>0$. The argument
is similar to that above. By Lemma 12 we get exponential coordinates on balls for the $\hat{g}$ metric on $X^{c}_{E}$ whose image in $\oP$ contains a fixed euclidean neighbourhood of $q$.
Arguing just as in the previous case, we get bounds on the covariant derivatives
of the metric tensors and can suppose that, in geodesic coordinates these
converge, along with the complex structures. So we have $(g^{\alpha}, J^{(\alpha)})\rightarrow
(g^{\infty)}, J^{\infty})$ say. For each $\alpha$ we have a pair of $J^{(\alpha)}$-holomorphic, commuting, vector fields $v_{1}^{(\alpha)}, v_{2}^{(\alpha)}$ and $Iv_{1}$
is a Killing field generating a circle action fixing $q$. Just as in the
previous case, the exponential map is equivariant for this action so the
limiting metric $g^{(\infty)}$ is also preserved by the same fixed circle
action. For the other sequence of vector fields $v_{2}^{(\alpha)}$ we have
to argue differently. We know that these are bounded in $L^{\infty}$ so it
follows from the ellipticity of the $\db$-operator that we can suppose
(after perhaps taking a subsequence) that these converge. What we have to
see is that the limit $v_{2}^{(\infty)}$ is not a multiple of $v_{1}$. But
this is the case, since
$v_{1}$ vanishes at $q$ while the length of $v_{2}(q)$ is $u^{22}$ which is
bounded below by Lemma 4. So we obtain, in the limit in geodesic coordinates over
a small neighbourhood of $q$
\begin{itemize}
\item a Kahler metric $g^{(\infty)}, J^{(\infty)}$;
\item a pair of linearly independent, commuting,  holomorphic vector fields $v_{1}^{(\infty)},
v_{2}^{(\infty)}$   such that $Iv_{1}^{(\infty)}$ generates the standard circle action.
\end{itemize}
Then just as before it follows from Guillemin's analysis that this data corresponds
to a function satisfying Guillemin boundary conditions on a neighbourhood
of $q$ in $\oP$.

\subsection{Rescaling}
It is standard practise in Riemannian geometry to rescale a metric in order
to obtain a fixed bound on the curvature.  We want to implement this idea in
our special situation. Suppose $u$ is a convex function on a polygon $P$ which
satisfies Guillemin boundary conditions defined by weights $\sigma$, with
$u^{ij}_{ij}=-A$. Let
$\lambda$ be a positive real number. Define a function $\tu$ on the polygon
$\tP= \lambda P$ by
$$  \tu(x^{1}, x^{2})= \lambda u(\lambda^{-1} x^{1}, \lambda^{-1} x^{2}).
$$
\begin{prop}
\begin{itemize}
\item The function $\tu$ satisfies Guillemin boundary conditions for the
weights $\tilde{\sigma}(\lambda E)= \lambda \sigma(E)$.
\item The curvature $\tilde{F}$ of $\tilde{u}$ satisfies
$$  \vert \tilde F \vert (\lambda p) = \lambda \vert F \vert(p). $$
\item The scalar curvature $\tA=\tu^{ij}_{ij}$ is
$$  \tA(\lambda p) = \lambda A(p). $$
\item If $u$ satisfies an $M$-condition then so does $\tu$ (with the same
value of $M$).
\item $\mu_{\tP, \tilde{\sigma}}= \mu_{P,\sigma}$.
\end{itemize}
\end{prop}

All of these are very easy to check. Notice that the second item implies
that
\begin{equation}  \int_{\tP} \vert \tF\vert^{2} d\mu_{Euc} = \int_{P} \vert F \vert^{2}
d\mu_{Euc}. \end{equation}
If $(P,\sigma)$ is Delzant, so also is $(\tP, \tilde{\sigma})$. There is then a
canonical diffeomorphism from $X^{c}(P,\sigma)$ to $X^{c}(\tP, \tilde{\sigma})$
and under this the Riemannian metric $\hat{g}$ is scaled by a factor $\lambda$.
In this case (5) is just the standard fact that the $L^{2}$ norm of the curvature
tensor is scale invariant in four real dimensions. 

Using this rescaling we can transfer the results of Section 3, under the
hypothesis that $\vert F \vert \leq 1$, to the general case. In fact  we have the following refinement of Proposition
4.
\begin{prop}
Let $(P,\sigma)$ be Delzant and let $\hat{g}$ be a metric on $X^{c}=X^{c}(P,\sigma)$
determined by a convex function $u$ on $P$. Suppose $u$ satisfies an $M$-condition. There is a $c>0$, depending only on $M$ and $\mu(P)$, with the following
property. For any $\rho>0$ and point $x\in X^{c}$

\

{\it either} there is a point $x'\in  X^{c}$ with $\distg(x,x')\leq \rho$  
and $\vert F(x') \vert \geq \rho^{-2}$,

\

{\it or} $\vert F(x')\vert \leq \rho^{-2}$ for all $x'$ with $\distg(x,x')\leq
\rho$ and  the exponential map at $x$ is an embedding on the ball of radius
$c\rho$ \end{prop}
This follows from Proposition 4 after rescaling and the observation that the hypothesis
$\vert F\vert\leq 1$ in Proposition 4 is only used on points within a fixed distance
of $x$.

\subsection{Blow-up limits}
Now suppose that, in our sequence $u^{(\alpha)}$ as considered in Theorem
1, the curvature $\vert
F\vert$ does not satisfy a uniform bound.
 For each $\alpha$ choose a point $p_{\alpha}$ where the modulus of the curvature achieves
its maximal value $K_{\alpha}$ and suppose that $K_{\alpha}\rightarrow \infty$.
We want ultimately to derive a contradiction. By translation we can suppose
that each $p_{\alpha}$ is the origin. We dilate by a factor $K_{\alpha}$
so we get a new sequence of data $(\tP^{(\alpha)}, \tsigma^{(\alpha)}$ and
functions $\tu^{(\alpha)}$.   It is clear that, perhaps after taking a subsequence,
one of three cases must occur. 
\begin{itemize}
\item The limit of the $\tilde{P}^{(\alpha)}$ is the whole of  $\bR^{2}$;
\item The limit of the $\tilde{P}^{(\alpha)}$ is a half-plane; 
\item The limit if the $\tilde{P}^{(\alpha)}$ is  a quarter-plane ({\it i.e.}
a nontrivial intersection of two half-planes).
\end{itemize}
(Here by the statement that \lq\lq the limit of $\tilde{P}^{(\alpha)}$ is $G$'' we mean that point 
of $G$ is contained in $\tilde{P}^{(\alpha)}$ for all large enough $\alpha$ and any point
not in the closure of $G$ is in the complement of $\tilde{P}^{(\alpha)}$ for all large enough $\alpha$.)

The main result of this subsection is
\begin{prop}
 If the limit of $\tilde{P}^{(\alpha)}$ is $G$, for one of the three cases above,
 then after taking a subsequence and adding suitable affine linear functions
 the $\tilde{u}^{(\alpha)}$ converge to a smooth convex function $\tilde{U}$ on
 $G$ which satisfies the equation $\tU^{ij}_{ij}=0$. The limit $\tU$ satisfies
 an $M$ condition in $G$. In the case  when $G$ is a quarter plane, the limit
 satisfies Guillemin boundary conditions and defines a complete, non-flat,
  zero scalar
 curvature  Kahler metric
 on $\bR^{4}$ with curvature in $L^{2}$.
 
 \end{prop}
 
 We give the proof in the three cases.
 
 \
 
 {\it Case 1:  The limiting domain is the whole plane.}
 
 \
 
 We can apply the results from Section 3 to the functions $\tu=\tu^{(\alpha)}$. We want
 to show that on any compact subset $K\subset \bR^{2}$ we have upper and
 lower bounds
 $$  C_{K}^{-1} \leq  \tu_{ij} \leq C_{K}. $$ 
 
 The upper bound follows immediately from Lemma 4 (since on compact sets
 the Euclidean distance to the boundary of $\tP^{(\alpha)}$ tends to
 infinity with $\alpha$). Let $J=J^{(\alpha)}$ be the function $\det(\tu_{ij})$.
 The crucial thing is to get a lower bound on $J(0)$. Corollary 2  implies that the distance in the metrics
 $\tilde{g}^{\alpha}$ corresponding to $\tu^{\alpha}$ from the origin to
 the boundary of $\tP^{\alpha}$ tends to infinity. By construction, $\vert
 \tilde F^{(\alpha)}\vert$ is equal to $1$ at the origin. We want to apply
 Proposition 5. Notice that when we rescale the derivatives of the scalar curvature function
 decrease, 
 so are certainly uniformly bounded in the sequence. Thus by
   Proposition 5 we can find a fixed small
 number $\delta$ such that $\vert F^{(\alpha)} \vert\geq 1/2$ on the $\tilde{g}^{\alpha}$
 ball of radius $\delta$ about the origin. On the other hand Lemma 8 implies that
 this ball contains a Euclidean ellipse of area at least $c J(0)^{-1/2} \delta^{2}$,
 for some fixed $c$. Thus 
 $$  \int_{\tilde P^{\alpha}} \vert \tilde{F}^\alpha \vert^{2} d\mu_{Euc}
 \geq c \delta^{2} J(0)^{-1/2}. $$
 Since, from (4) and (5),  the integral on the left is bounded, we obtain a lower
 bound on $J(0)$, as required. Combined with the upper bound on $\tu_{ij}$ this
 lower bound on $J(0)$ yields an upper on $\tilde{u}^{ij}$ at the origin. Now
 Lemma 7 gives an upper bound on $\tu^{ij}$ at points of bounded $\tilde{g}$ distance
 from the origin. The upper bound on $\tu_{ij}$ implies that on compact subsets
 of the plane the $\tilde{g}$ distance to the origin is bounded. So we conclude
 that $\tu^{ij}$ is bounded above on compact subsets of the plane, which is
 the same as the lower bound on $\tu_{ij}$.
 Once we have these upper and lower bounds on $\tu_{ij}$ the convergence of a subsequence is straightforward, just
 as in the proof of Proposition 8, and the fact that the limit $\tU$ has $\tU^{ij}_{ij}=0$ follows
 from the third item of Proposition 9.

\

{\it Case 2: The limiting domain is a half-plane.}
 
 \
 
 The proof is similar to the first case. The upper bound on $u_{ij}$ on compact
 subsets of the limiting half-plane is obtained just as before. Let $d_{\alpha}$
 be the distance from the origin to the boundary of $\tP^{(\alpha)}$ if $d_{\alpha}$
 is bounded below we can argue just as before. The only difficulty comes
 when $d_{\alpha}\rightarrow 0$, which is the same as saying that the origin
 is on the boundary of the limiting half-plane. Fix a parameter
 $\tau<\pi/2$. For each $\alpha$ we take
 a point $q_{\alpha}$ on the boundary of $\oP^{\alpha}$ which minimises the
 $g_{\alpha}$ distance to the origin and let $p_{\alpha}$ be the point of
 $\tP^{\alpha}$ a distance $\tau$ from $q_{\alpha}$ along the geodesic emanating
 from $q_{\alpha}$ orthogonal to the boundary of $\tP^{\alpha}$. Then by
 Lemma 13 the distance from $p_{\alpha}$ to the boundary of $\tP^{\alpha}$ is at
 least $\tau$ (once $\alpha$ is sufficiently large). Here we use the fact
 that the distance from the origin to all but one of the edges of $\tP^{\alpha}$
 tends to infinity with $\alpha$. Now by applying Proposition 5 to a geodesic ball centred
 at $q_{\alpha}$ we see that we can fix $\tau$ so that $\vert F\vert \geq
 1/2$, say, on the ball of radius $\tau/2$ about $p_{\alpha}$. Now the argument
 goes through just as before.    

\

{\it  Case 3: The limiting domain is a quarter-plane.}

\

 The proof in this case is much like that of Proposition 8. Let $q_{\alpha}$ be the vertex
 of $\tP^{\alpha}$ close to the vertex of the limiting quarter-plane and
 let $E_{1}^{\alpha}, E_{2}^{\alpha} $ be the edges of $\tP^{\alpha}$ meeting
 in $q_{\alpha}$ with defining functions $\lambda_{i,\alpha}=\lambda_{E_{i}^{\alpha}}$.
 Observe that the definition of the $\tsigma^{\alpha}$ implies that that these
 $\lambda_{i,\alpha}$ converge as $\alpha$ tends to infinity to defining functions
 for the edges of the quarter plane.  We
 obtain the lower bounds on $\tu_{ij}$, or equivalently the upper bound on
 $\tu^{ij}$, by applying Lemma 6, using  the geodesics emanating from $q_{\alpha}$,
 and the upper bounds on $u_{ij}$ using the $M$-condition and Lemma 4. Just as
 in the proof of Proposition 8 we show that the limit satisfies Guillemin boundary conditions along the edges of
 the quarter plane, and it is clear that the corresponding $4$-manifold is
 diffeomorphic to $\bR^{4}$. The completeness of the limiting metric follows
 from general principles or more directly from our estimate
 $$ \lambda_{i,\alpha}\leq \cosh(\distg(q_{\alpha}, p))-1 . $$
 The fact that the curvature of the limiting metric is in $L^{2}$ follows
 from (4), (5) and Fatou's Lemma. Of course, the fact that the limiting metric is not flat
 follows from the normalisation that $\vert \tF^{\alpha}\vert$ is equal to $1$ at the
 origin, and the $C^{\infty}$ convergence.

 With Proposition 11  in place the desired contradiction (to the hypothetical
  blow up of the curvature
 in the sequence) follows from the following two results.
 
 \begin{thm}
 There is no convex function $U$ on a half-plane which satisfies an $M$-condition
 and the equation $U_{ij}^{ij}=0$.
 \end{thm}

 \begin{thm}
 If $U$ is a convex function on a quarter plane which satisfies an $M$-condition
 and  which defines a complete
 zero scalar curvature metric on $\bR^{4}$ with curvature in $L^{2}$ then
 the metric is flat.
 \end{thm}
 
 We will give one proof of Theorem 3 now. This uses a result of Anderson \cite{kn:An},
 which we quote.
 \begin{thm}[Anderson] Let $g$ be a  complete self-dual Riemannian metric on $\bR^{4}$ with zero
 scalar curvature. Suppose that the curvature of $g$ is in $L^{2}$ and that
 the volume $V(r)$ of the ball (in the metric $g$) of radius $r$ about the
 origin satisfies $V(r)\geq c r^{4}$ for some $c>0$. Then $g$ is flat.
 \end{thm}
 To see that this applies to our case, recall first that scalar-flat Kahler
 metrics in two complex dimensions are self-dual. Thus the only thing we
 need to establish is the volume growth. This uses the $M$-condition. We
 can suppose the quarter plane in question is the standard one defined by
 $x^{i}\geq 0$ and that the boundary conditions correspond to the defining
 functions $x^{i}$. For $\tau>0$ let $\Omega(\tau)$ be the triangle
 $\{ (x^{1}, x^{2}) : x^{i} \geq 0, x^{1}+x^{2} \leq \tau\}$ and let $X^{c}(\tau)$
 be the corresponding subset of $X^{c}$. By Corollary 1 we have $\distg(p,0)\leq C
 \sqrt{\tau}$ for $p\in \Omega(\tau)$, where $C$ depends on $M$. So $X^{c}(\tau)$
 is contained in the ball of radius $C \tau^{2}$. On the other hand the volume
 of $X^{c}(\tau)$ is equal to $(2\pi)^{2}$ times the Euclidean area of $\Omega(\tau)$
 which is $\tau^{2}/2$. So we deduce that
 $$  V(C \sqrt{\tau}) \geq 2 \pi^{2} \tau^{2}, $$
 from which the statement follows.
\section{Nonexistence of blow-up limits}
\subsection{The case of the half-plane}

Throughout this section we will, contrary to our general convention, use
lower indices for our coordinates $(x_{1}, x_{2})$ on the Euclidean plane.

We will first indicate the proof of Theorem 2 in the case of a convex function $u$ satisfying
the zero scalar curvature equation
$u^{ij}_{ij}=0$ and an $M$-condition on the whole plane. While this is subsumed
in the harder case below the proof is substantially simpler. We suppose $u$
is normalised to achieve its minimum at the origin. Then an easy elementary
argument (see the proof of Lemma 14 below) shows that an  $M$ condition implies a uniform bound
on the first derivative, $\vert \nabla u \vert \leq M'$, say, on $\bR^{2}$.
Now we apply Theorem 5 of the Appendix to the restriction of $u$ to a large Euclidean
disc of radius $R$ centred at the origin. This yields $\det (u_{ij})(0)\leq
CR^{-2}$, and we get a contradiction by letting $R$ tend to infinity.

Now we give the proof for the case when the function is only defined on a
half-plane.  
\begin{lem}
Suppose $(p_{1}, p_{2})$ is fixed and $u$ is a convex function on a neighbourhood of a rectangle
$\{ x_{1}, x_{2}:\vert x_{1}- p_{1}\vert \leq L_{1}, \vert x_{2}-p_{2}\vert \leq L_{2}\}$
which satisfies the zero scalar curvature equation $u^{ij}_{ij}=0$. Let
$$ V_{1}= V((-L_{1},0),( L_{1}, 0)), V_{2}= V((0,-L_{2}), (0,L_{2})) $$
and set $\Delta= {\rm Max} (V_{1} L_{1}, V_{2} L_{2})$. Write $J$ for the
function $\det(u_{ij})$.
Then there is a universal constant $\kappa$ such that
$$   J(0,t)\leq  \frac{\kappa \Delta^{2}}{L_{1}^{2} L_{2}^{2}}, $$
for $\vert t\vert \leq L_{1}/4$,
\end{lem}
(Recall that the function $V(p,q)$ is defined in  (2) in Section 1.)

Obviously we can suppose $(p_{1}, p_{2})=(0,0)$.
It is elementary to check that the statement is invariant under dilations
of the co-ordinates, so we can reduce to the case when $L_{1}=L_{2}=1$. We
can also suppose that $u$ is normalised so that it vanishes, together with
its first derivatives, at the origin. So $u$ is positive and, by convexity
and the definition of $V_{1}$,
the modulus of the partial derivative $\frac{\partial u}{\partial x_{1}}$ is bounded by
$\Delta$ on the interval $\{(t,0): -1\leq t\leq 1 $. Similarly for the $x_{2}$
variable. Thus $u\leq \Delta$
at the four points $(\pm 1,0), (0,\pm 1)$. By convexity, $u\leq \Delta$ on
the square $K$ formed by the convex hull  of these four points.
Let $D$ be the disc of radius $\thalf$ about the origin. So $D$ is contained
in the interior of $K$ and the distance from $D$ to the boundary of $K$ is
$d=(1/\sqrt{2})-\thalf$. By an elementary property of convex functions we have
$\vert \nabla u\vert \leq \rho =\Delta/d $ on $D$. Then by Theorem 5 of the
Appendix there is a universal constant $C$ such that $J\leq C \rho^{2}$ on
the interior disc of radius $\tquart$ centred on the origin. Thus we can take
$\kappa= Cd^{-2}$.

\begin{lem}
Suppose $u$ is a convex function on the half-plane $\{ (x_{1}, x_{2}): x_{1}>0\}$
which satisfies the $M$ condition and the zero scalar curvature equation
$u^{ij}_{ij}=0$. Write $J=\det (u_{ij})$. Then
\begin{enumerate}
\item For any $\eta>0$ there is an $h_{0}$ such that $J(x_{1}, x_{2}) \leq
\eta x_{1}^{-1} $ if $x_{1} \geq h_{0}$.
\item For any $\epsilon>0$ and $0< h_{1} < h_{2}$ there is a $T>0$ such that
$J(x_{1}, x_{2}) \leq \epsilon $ if $\vert x_{2}\vert \geq T$ and $h_{1}\leq x_{1}\leq
h_{2}.$
\end{enumerate}
\end{lem}

To prove the   first item we consider a point $p=(p_{1},p_{2})$ with $p^{1}>0$.
and consider the rectangle $\{(x_{1}, x_{2}): \vert x_{1}- p_{1}\vert
\leq p_{1}/2, \vert x_{2}-p_{2} \vert\leq p_{1}\}$. We can apply Lemma 14, where
$L_{1}= p^{1}/2, L_{2}= p^{1}$. The $M$ condition implies that $V_{1}\leq 2M, V_{2}\leq
M$, so $\Delta \leq M p_{1}$. Then we obtain $J(p) \leq 4 \kappa M^{2} (p_{1})^{-2}$
and the result follows (with $h_{0}=4 \kappa M^{2} \eta^{-1}$).

To prove the second item we first consider the case when $h_{1}= \frac{3h}{4},
h_{2}= \frac{5h}{4}$ for some $h$. It obviously suffices to show that the
statement is true for $x^{2} \geq T$, once $T$ is suitable large. The $M$
condition implies that 
$$  \int_{-\infty}^{\infty} u_{22}(h, t)  dt \leq M< \infty,  $$
so given any $\delta$ we can find a large $S$ such that 
\begin{equation}  \int_{S}^{\infty} u_{22}(h,t) dt \leq \delta. \end{equation}
Now, for $L_{2}>0$ consider  a  rectangle
$$  \{ (x^{1}, x^{2}): \frac{h}{2}\leq x_{1}\leq \frac{3h}{2}, \vert x_{2}-p_{2}\vert \leq L_{2}\}
$$ as before. Suppose that $p_{2}-L_{2}\geq S$. Then (6) implies that $V_{2}\leq \delta$. Suppose that $\delta L_{2}\leq V_{1} L_{1}= M h $. Then $\Delta=M
h$ and Lemma 14 gives $J( x_{1}, p_{2}) \leq 4 \kappa^{2} M^{2} / L_{2}^{2}$ for
$3h/4\leq x_{1} \leq 5 h/4$.  So,
given $\epsilon>0$ we first choose a large $L_{2}$ such that
$4\kappa^{2} M^{2}/L_{2}^{2} \leq \epsilon$. Then we choose a small $\delta$
such that $\delta L_{2} \leq M h$. Then we choose $S$ as above and set $T=
S+ L_{2}$. 

Finally, for general $h_{1}< h_{2}$, we cover the interval $[h_{1}, h_{2}]$
with a finite number of intervals of the form above and take the maximum
value of the corresponding $T$'s.

Now we can prove Theorem 2. Suppose that $u$ satisfies the $M$-condition on the half-plane
and $u^{ij}_{ij}=0$. We consider the function $F=\det(u_{ij})^{-1}$ on the
half-plane. This satisfies the equation
$$  u^{ij}F_{ij}=0. $$
(See \cite{kn:D2}, Sec. 2.1). So, for any constant $\lambda$, the function $G=F-\lambda x_{1}$ satisfies
$u^{ij} G_{ij}=0$, hence can have no local minumum or maximum.  Suppose,
without loss of generality, that $F(1,0)= 1$ and take $\lambda=2$, so
$G(1,0)=-1$. Clearly $G(x_{1}, x_{2})\geq -\thalf$ if $x_{1}\leq \tquart$. By the first
item in Lemma 15 we can choose $h_{0}$ so large that $F(x_{1}, x_{2}) \geq
10 x_{1}$ if $x_{1}\geq h_{0}$. This implies that $G(x_{1}, x_{2})\geq 8
x_{1}>0$ if $x_{1}\geq h_{0}$. By the second item of Lemma 15 we can choose
$T$ so large that $F(x_{1}, x_{2}) \geq 10 h_{0}$ if $\vert x_{2} \vert \geq
 T$ and $\frac{1}{4} \leq x^{1}\leq h_{0}$. This implies that $G(x^{1}, x^{2})
 \geq 8 h_{0} >0$ if $\vert x^{2}\vert \geq T$ and  $\frac{1}{4} \leq x^{1}\leq h_{0}$. So $G(x^{1}, x^{2})
 \geq -\thalf$ if $(x_{1}, x_{2}) $ lies on the boundary of the rectangle
 $$ Q= \{(x_{1}, x_{2}): \vert x_{2} \vert \leq T, \tquart\leq x_{1}\leq h_{0}
 \}. $$
 Since $G$ takes the value $-1$ on the interior point $(1,0)$ of $Q$ it must
 have an interior minumum, which is the desired contradiction.
 
 \

 One point worth noting here is that in the case of a function defined on
 the whole plane the argument can be made entirely effective. There is no
 need to take the limit as $\alpha$ tends to infinity of the sequence $\tu^{(\alpha)}$
 in Proposition 11.  The same argument
 can be used to obtain an explicit {\it a priori} estimate of the form
 $$   \vert F \vert \leq C \distg(\ , \partial P)^{-2}. $$
 In the case of the half-plane it seems to be necessary to pass to the limit,
 and the proof does not yield an explicit {\it a priori} estimate in general.
 However if one considers the case when $A^{(\alpha)}\geq 0$ then this can
 be done, and one gets an  explicit estimate of the form
 $$ \vert F \vert \leq C \distg(\ , V)^{-2}, $$
 where $V$ is the set of vertices.

\subsection{The case of the quarter-plane}
Here we give a second, self-contained, proof of Theorem 3.
The general strategy of the proof is in part similar to Anderson's,  in that
we show that the curvature tensor vanishes by applying an integral formula
for its $L^{2}$ norm,
and the crux of the matter is to establish that the relevant boundary term vanishes
in the limit. We will first state the relevant integral formula, in our special
situation.

Let $u$ be a convex function on $\bR^{+}\times \bR^{+}$ with $u^{ij}_{ij}=0$,
as in the statement
of the Theorem, and for $R>0$ let $\Omega(R)$ be the triangle formed by the
intersection of the quarter plane with the half-space $\{x_{1}+x_{2} \leq
R\}$. Let $\partial \Omega(R)$ denote the ordinary boundary of the triangle,
made up of three line segments and $\partial^{0} \Omega(R)$ be the single
segment lying on the line $\{x_{1}+x_{2}=R\}$. (Recall that $\Omega(R)$ corresponds
to a differentiably embedded ball in the $4$-manifold $X^{c}$ and $\partial^{0}
\Omega(R)$ corresponds to the boundary of this ball. The other two segments
in $\partial \Omega(R)$ correspond to fixed points for the two basic circle
actions on $X^{c}$.) Now we have

\begin{equation} \int_{\Omega(R)} \vert F\vert ^{2} d\mu = \int_{\partial^{0} \Omega(R}  \nu^{i}
\end{equation} where 
\begin{equation} \nu^{i}= F^{ij}_{ab} u^{ab}_{i} - F^{ia}_{ib} u^{bj}_{j} .\end{equation}
 The integrand on the right hand side of the formula (7) is
written as a vector field but this can be viewed as a $1$-form using the
canonical identification furnished by the Euclidean area element $d\mu= dx_{1}
dx_{2}$. We leave the verification  of this identity as an exercise for the
reader (see also the similar discussion in \cite{kn:D2}, Sec. 5.2). The overall strategy of our
proof is to show that the integral on the right hand side of (7) tends to zero
as $R\rightarrow \infty$, which implies that $F$ is identically zero. 

We begin by establishing that the curvature decays as a function of the Riemannian
distance from the origin. To fit in with the wider literature we will phrase
this discussion in terms of the Riemannian $4$-manifold $X^{c}$, although
of course it can be translated into the two-dimensional language. The crucial
thing is that this Riemannian manifold has the  property stated in Proposition
10 (The discussion there assumed a compact manifold but it is
easy to see that the proofs work equally well in the present situation.)
Moreover if (in the notation of Proposition 10) $\vert F\vert \leq \rho^{-2}$
on the ball of radius $\rho$ about a point $x$ in $X^{c}$ we have, by applying
Proposition 5,
$$ \vert F(x) \vert^{2} \leq C \rho^{-4} \int_{B(c\rho,x)} \vert
F\vert^{2} dV $$ for some fixed $C$, where $dV$ is the Riemannian
volume element. Now we recall a general fact:
\begin{lem} Let $X$ be a complete, noncompact, Riemannian manifold with
base point $x_{0}$. Let   $K$ be a continuous, non-negative, $L^{2}$ function on $X$
with the following property. There  are constants $c,C>0$ such that for any $x\in X$ and $\rho>0$, either there is
a point $x'\in X$ with $d(x,x')\leq \rho$ and $K(x')> \rho^{-2}$ or 
$$  K(x)^{2}\leq C \rho^{-4} \int_{B(c\rho, x)} \vert K \vert^{2}. $$
Then   
$K(x) d(x,x_{0})^{2}\rightarrow 0$ as $x$ tends to infinity in $X$.
\end{lem}

To see this, let $E>0$ and let $X_{E}\subset X$ be a compact set such that
$$  \int_{X\setminus X_{E}} K^{2} \leq C^{-1} E. $$
Define a function $\rho_{E}$, taking values in $(0,\infty]$, by
$$ \rho_{E}(x)^{-4}=  \frac{K(x)^{2}}{ 2 E}. $$
It is convenient to work with this and one can check step-by-step in the
argument below that, with  the obvious interpretations, there are no problems
from the zeros of $K$. The crucial thing is that $\rho_{E}$ is bounded below
by a strictly positive number on any compact set in $X$. Now suppose $x$
is a point in $X$ with $\rho_{E}(x)\leq \epsilon d(x,X_{E})$ where $\epsilon<c$.
Then the ball of radius $c \rho_{E}(x)$ about $X$ does not meet $X_{E}$ so
the second alternative in the hypothesis (taking $\rho=\rho_{E}(x)$) would give $K(x)^{2} \leq \frac{1}{2}
K(x)^{2}$. Since $\rho_{E}(x)$ is finite $K(x)$ is nonzero and we conclude
that the first alternative must hold; that is, there is a point $x'$ with
$$ d(x,x') \leq \rho_{E}\ \ , \ \ K(x')>\rho_{E}(x)^{-2} = K(x)/\sqrt{2E}.
$$
So now we have $\rho_{E}(x')\leq 2^{-1/4} \rho_{E}(x)$. We also have
$$\epsilon d(x',X_{E}) \geq \epsilon d(x, Z_{E})- \epsilon d(x,x')\geq (1-\epsilon)
\rho_{E}(x). $$
Thus $$\rho_{E}(x')\leq \frac{\epsilon}{(1-\epsilon)2^{1/4}} d(x',Z_{E}).
$$
Suppose $\epsilon$ is so small   that $2^{1/4}(1-\epsilon) >1$. Then
$\rho_{E}(x') \leq \epsilon d(x',X_{E})$.
Thus $x'$ satisfies the same hypothesis as $x$ did. We continue in this way
to generate a sequence $x_{n}$ with 
$$  d(x_{n}, x_{n+1}) \leq \rho_{E}(x_{n}), $$
and $$\rho_{E}(x_{n+1} \leq  2^{-1/4}\rho_{E}(x_{n}). $$
Thus $x_{n}$ is a Cauchy sequence in $X$ and $\rho_{E}(x_{n})$tends to zero,
a contradiction. So we conclude that for $\epsilon< {\rm Min}(c,1-2^{-1/4})$
we have $$\rho_{E}(x)\geq \epsilon^{-1} d(x,X_{E})$$ for all $x$ in $X$.
This says that
$$    K(x) d(x,X_{E})^{2} \leq 2 E/\epsilon, $$
and the result follows, since we can take $E$ as small as we please.

\

So in our case we know that the function $\vert F\vert$ on the quarter-plane
decays faster than than the inverse square of the Riemannian distance to
the origin. The next step is to relate this distance to the Euclidean distance
in the quarter-plane. For this we use another integral identity. Change Euclidean
coordinates by setting $ y=x_{1}-x_{2}, z=x_{1}+x_{2}$ and denote derivatives
with respect to the new coordinates by $u_{yy}$ {\it etc.} 
\begin{lem}
If $u$ satisfies $u^{ij}_{ij}=0$ in the quarter plane and Guillemin boundary
conditions then for any $R>0$
$$ \int_{y=-R}^{R} u^{zz}(y,R) dy = R^{2}. $$
\end{lem}

(The notation is slightly ambiguous here, so we should emphasise that in
the formula above we are regarding $u=u(y,z)$ as a function of $y,z$. The
region of integration is exactly  $\partial^{0}\Omega_{R}$, as  considered
above.)

To see this, let $f(y,z)$ be the function $f(y,z)= R-y$ on the region $\Omega_{R}$.
The zero scalar curvature condition takes the same form in the new coordinates,
so we write it as $u^{\alpha \beta}_{\alpha \beta}=0$, where $\alpha, \beta$
run over the labels $y,z$. So we have
$$ \int_{\Omega_{R}} u^{\alpha \beta}_{\alpha \beta} f = 0. $$
Now we integrate by parts twice. Since $f$ is linear we have $f_{\alpha \beta}=0$
and there is no contribution from the interior so we get the identity
$$ \int_{\partial \Omega_{R}} u^{\alpha \beta}_{\alpha} f = \int_{\partial
\Omega_{R}} u^{\alpha \beta} f_{\alpha}. $$  The function $f$ vanishes on
$\partial^{0} \Omega(R)$ and the Guillemin boundary conditions imply that $u^{\alpha \beta}_{\alpha}$
has normal component $1$ along the axes. Thus
$$  \int_{\partial \Omega(R)} u^{\alpha \beta}_{\alpha} f = 2\int_{0}^{R}\
R-t \ dt= R^{2}.$$
On the other hand the boundary conditions imply that the normal component
of $u^{\alpha \beta} f_{\alpha}$ vanishes along the axes, so
$$ \int_{\partial \Omega(R)} u^{\alpha \beta} f_{\alpha} = \int_{\partial^{0}
\Omega(R)} u^{\alpha \beta} f_{\alpha},$$ and the result follows.

Now for any fixed $z>0$ let $s(z)$ be the Riemannian distance from the origin
to the interval $\partial^{0}\Omega(z)= \{x^{1}+x^{2}=z\}$. Suppose for the moment that this distance
is realised by a unique  minimal geodesic $\gamma$ and that there is no Jacobi
field along $\gamma$ which vanishes at the origin and is tangent to this
interval at the other end point. Then $s$ is smooth around this value of
$z$ and 
$$    \frac{ds}{dz} =  \vert dz \vert_{g}^{-1} = \left(u^{zz}\right)^{-1/2}, $$
where $u^{zz}$ is evaluated at the distance-minimising point of the interval.
In any case, $s(z)$ is a Lipschitz function and if we define 
$$   \phi(z)= \max u^{zz}, $$
where the maximum is taken over this interval, then we have
\begin{equation} \frac{ds}{dz} \geq \phi(z)^{-1/2},  \end{equation}
interpreted in an appropriate generalised sense. 

We want to go from the integral identity of Lemma 17 to a pointwise bound on $u^{zz}$,
and hence on $\phi(z)$.
For this we use
\begin{lem} 
Suppose $f$ and $\sigma$  are positive function on an interval $[0,R]$, where
 $R\geq 1$, with
$$ \vert  f''(t) \vert \leq f(t) \sigma(t), $$
and for any $\lambda>0$ we have
$$ \int_{\lambda/2}^{\lambda} \sigma(t) dt \leq 1. $$
Then for any $t_{0}$ in $[0,R]$ we have $$f(t_{0}) \leq 18 \int_{0}^{R} f(t) dt.
$$
\end{lem}
To simplify notation we will give the proof in the case when $f$ attains
its maximum at $t_{0}=0$. It
will be clear that this is the \lq\lq worst'' case and that the argument
applies to all points. Set
$$ I=\int_{0}^{R}  f(t) dt. $$

For $h>0$ use the formula
$$  f(0)= f(h)- h f'(h) - \int_{0}^{h}t f''(t) dt. $$

Using the assumption that $f$ attains its maximum at $0$, and the given differential
inequality $\vert f''\vert\leq f \sigma$,  we have
$$\vert \int_{h/2}^{h} t f''(t) dt \vert \leq \frac{f(0) h}{2} \int_{h/2}^{h} \sigma (t) dt \leq\frac{ f(0) h}{2}. $$

Summing over a geometric series, as in the proof of Lemma 2,  we obtain
$$ \vert \int_{0}^{h} t f''(t) dt \leq h f(0). $$
So  if $h\leq \thalf$, say, we have
\begin{equation} f(0)\leq 2( f(h)- hf'(h)). \end{equation}
Let $t_{1}>0$ be a point in the interval $[0,\tthird]$ where $f$ attains its minimum.
Then $$  I\geq \int_{0}^{1/3} f(t) dt \geq \frac{f(t_{1})}{3}. $$
If $t_{1}<\tthird$ the derivative $f'(t_{1})$ vanishes and, taking $h=t_{1}$
in the inequality above we have
$$  f(0) \leq 2 f(t_{1}) \leq 6 I. $$
Suppose, on the other hand, that $t_{1}=\tthird$. Let $g$ be the affine-linear
function with $g(\tthird)=f(\tthird)$ and $g(\thalf)=f(\thalf)$. If $f'(\tthird)<g'(\tthird)$ then
there is a point $h$ in the interval $(\tthird,\thalf)$ where $f'(h)=g'(h)$ and
$f(h)<g(h)$. If $f'(\tthird)\geq g'(\tthird)$ we take $h=\tthird$. In either case
$$f(h)- hf'(h)\leq g(0)=  \frac{ \frac{1}{2} f(\tthird)- \frac{1}{3} f(\thalf)}{\frac{1}{2}-\frac{1}{3}}
\leq 3 f(\tthird).$$ 
Then applying the inequality (10) above with this value of $h$ we obtain
  $$ f(0) \leq 6f(1/3)\leq 18 I. $$

  \begin{cor}
  Suppose that $u$ satisfies $u^{ij}_{ij}=0$ in $\bR^{+}\times \bR^{+}$ and Guillemin
  boundary conditions. Suppose that $u$ satisfies the $M$ condition with $M=1$
  and that, for some $R>1$, $\vert F\vert\leq 1$ on $\partial\Omega_{R}$. Then
  $$  u^{zz}\leq 18 R^{2} $$
  on $\partial^{0} \Omega_{R}$.
  \end{cor}
  To see this, observe that
  $$  \vert u^{zz}_{yy}\vert  \leq \vert F\vert u^{zz} \ u_{yy}. $$
  Then the result follows from Lemmas (17) and (18), taking $f(t)=u^{zz}(t+R,R)$
  and $\sigma(t)=u_{yy}(t+R,R)$.

 \begin{lem}
 There is a constant $C$ such that $z\leq C s(z)^{2}$ for all $z$.
 \end{lem} 
  
  For $z>0$ define $\lambda(z)= \frac{z}{s(z)^{2}}$. We know, by Corollary
  2, that $s(z)\rightarrow \infty$
  as $z\rightarrow \infty$. Further, we know that $\vert F\vert = o(s^{-2})$,
  so it follows that for any $\epsilon>0$ we can find an $R_{0}$ such that
  $$\vert F\vert \leq \epsilon \frac{\lambda(R)}{R}, $$
 on $\partial^{0} \Omega_{R}$,  once $R\geq R_{0}$. For a fixed $R\geq R_{0}$,
 suppose that $\lambda(R) \geq \epsilon^{-1}$ and 
 set 
 $$\eta= \frac{\epsilon \lambda(R)}{R}, $$
 so by hypothesis $\eta \geq R^{-1}$.
 Now rescale using this factor $\eta$, so we define $\tilde{u}(x_{1}, x_{2}) = \eta u(x-{1}/\eta, x-{2}/\eta)$.
 Set $\tilde{R}= \eta R$ and consider the rescaled solution on the triangle
 $\Omega(\tilde{R})$, which corresponds to the original solution on the triangle
 $\Omega(R)$.  
 The curvature tensor $\tilde{F}$ of $\tilde{u}$ satisfies $\vert \tilde{F}
 \vert \leq \eta^{-1} \vert F \vert \leq 1$ on $\partial^{0} \Omega_{R'}$ and
 so we can apply Corollary 3 to $\tilde{u}$ to get
 $$ \tilde{u}^{zz} \leq 18 \eta^{2} R^{2} $$
 on $\partial^{0} \Omega_{\tilde{R}}$. Transforming back, this becomes
 $$  u^{zz}\leq 18 \eta^{-1} \eta^{2} R^{2} = 18 \epsilon \lambda(R) R $$
 on $\partial^{0} \Omega_{R}$. In other words we have the following: for $R\geq R_{0}$ if $\lambda(R)
 \geq \epsilon^{-1} $ then $\phi(R) \geq (18 \epsilon \lambda(R) R)^{-1/2}$.
 Now consider the derivative of $\lambda$. Using (9) we have
 $$ \frac{d\lambda}{dz} = s(z)^{-2} \left( 1- \frac{2 z s'(z)}{s}\right)\leq
 s(z)^{-2}\left( 1- \frac{2 z \phi(z)}{s(z)}\right).$$
 If $z\geq R_{0}$ and $\lambda(z)\geq \epsilon^{-1}$ then we have
 $$ \frac{d\lambda}{dz} \leq s(z)^{-2}\left( 1- \frac{2 z}{s \sqrt{\epsilon
 \lambda z}}\right)=s(z)^{-2}(1- \frac{2}{\sqrt{18\epsilon}}. $$
 Now we fix $\epsilon<2/9= 4/18$ so that $(1- 2/{\sqrt(18\epsilon})<0$ and
 we see  that once $z\geq R_{0}$ and $\lambda>
 \epsilon^{-1}$ the function $\lambda$ is decreasing. It follows then that
 $\lambda(z)$ is bounded.

 Combining  Lemma 16  and Lemma 19, we have
 \begin{equation} \vert F \vert = o( z^{-1}) . \end{equation}
 For $R>1$ we now  rescale by $R$, so we define $\tilde{u}^{(R)}$
 (which we sometimes just denote by $\tilde{u})$) to be
  $$\tilde{u}^{(R)}(x_{1}, x_{2}) = R^{-1} u(R x_{1}, R x_{2})
 + L(x_{1}, x_{2}),$$
 where $L$ is an affine-linear function chosen so that $\tilde{u}^{(R)}$ and its
 first derivatives vanish at the point $x_{1}=\thalf, x_{2}=\thalf$. We consider
 the restriction of $\tilde{u}^{(R)}$ to the fixed quadrilateral
 $$ Q= \{ (x_{1}, x_{2}): \thalf< x_{1} + x_{2} < 2, \ \ x_{1}, x_{2}>0\}. $$
We write $\tilde{F}^{(R)}$ for the curvature tensor corresponding to $\tilde{u}^{(R)}$.
 The decay condition (11) implies that  $\vert \tilde{F}^{(R)}\vert $ 
 tends to zero on $Q$, as $R\rightarrow \infty$. As usual, we obtain an upper bound on the Hessian
 $\tu_{ij}$ over compact subsets of the interior of $Q$. Now Corollary 3 gives an upper bound on $\tilde{u}^{zz}$ over
 $Q$.  Lemma 3 gives
 $$   \frac{\partial^{2}}{\partial y^{2}} \left( \tilde{u}_{yy}^{-1}\right) \leq 1,
 $$
 say over $Q$.The boundary conditions fix the values of $\frac{\partial
 }{\partial y} \left (\tu_{yy}^{-1}\right) $ on $y=\pm z$ and this gives lower
  bound on $\tu_{yy}$,
 $$  \tu_{yy}^{-1}\leq C (z- \vert y \vert) $$
 Now
 $$  \tu^{zz} = \tu_{yy}/\det(\tu_{\alpha \beta}), $$
 so we obtain a lower bound on the determinant
 $$\det(\tu_{\alpha \beta})\geq C^{-1} (z-\vert y\vert)^{-1}. $$
  Combining with our upper bounds on the components of $u_{\alpha \beta}$
  we obtain upper bounds 
  $   \tu^{yy}\leq C, \vert \tu^{yz}\vert \leq C$ on $Q$.  Then, just as before, 
  we can conclude that as $R\rightarrow \infty$ the $\tu^{(R)}$ converge on compact
  subsets of the interior of $Q$ to a smooth
  limit $\tu^{(\infty)}$ with $\tilde{F}^{(\infty)}=0$. Now the boundary term
  in (7) is scale invariant, so we get the same computing with $\tu^{(R)}$
  and integrating over the fixed interval $\partial^{0}\Omega(1)$ in the
  interior of $Q$. It is then straightforward
  to check that this tends to  zero with $R$.

 \section{Appendix: applications of the maximum principle }
 
 In this appendix we use the maximum principle to derive upper and lower
 bounds on the determinant of the Hessian of a solution to Abreu's  equation.
 The results and their proofs are similar to those in \cite{kn:D2}, Sect.
 4,  but differ in being specific to the two-dimensional case. The inspiration
 for these results comes from the work of Trudinger and Wang in 
 \cite{kn:TW1},\cite{kn:TW2} and, particularly \cite{kn:TW3} , Remark 4.1.
 
 \begin{thm}
 Suppose that $u$ is a convex function on the closed disc of radius $R$ in
 $\bR^{2}$, smooth up to the boundary and with $\nabla u=0$ at the origin.
  Let $A(x)$ be the function $A=- (\partial u^{ij})_{ij}$ and
let $$A^{+}= \max (\max_{x}(A(x),0)), A^{-}= - \min(\min_{x} A(x), 0). $$
\begin{itemize}
\item If the derivative $\nabla u$ maps the $R$-disc to the disc $\vert \xi
\vert \leq \rho$ then on the interior disc $\{\vert x\vert\leq \tquart\}$ we
have
$$  \det (u_{ij}) \leq \left(\frac{\rho}{R}\right)^{2} \left( c_{1} + c_{2}
R^{2} \rho^{2} (A^{-})^{2}\right); $$
\item If the derivative $\nabla u $ maps the $R$-disc {\it onto} the disc
$\vert \xi \vert \leq \rho$ then on the set where $\vert \nabla u \vert \leq
\tquart$ we have
$$ \det(u_{ij})  \geq \left( \frac{\rho}{R}\right)^{2} \left( c_{3} + c_{4}
R^{2} \rho^{2}( A^{+})\right)^{-2}.;$$
for universal constants $c_{1}, c_{2}, c_{3}, c_{4}$. 
\end{itemize}
\end{thm}

Rescaling the domain and multiplying $u$ by a constant, we can assume that
$R=\rho=1$. We begin with the first item. Here we consider the function
$$ f= -L + F - \alpha g^{ab} u_{a} u_{b},$$
on the open disc, where $L= \log \det (u_{ij})$, $F$ is a smooth function
which tends to $+\infty$  on the boundary of the disc, to be specified shortly,
$\alpha$ is an arbitrary strictly positive constant and $g^{ij}$ denotes
the standard Euclidean metric tensor. (Thus $g^{ab} u_{a} u_{b}$ is another
notation for $\vert \nabla u\vert^{2}$.) The function $f$ attains its minumum
in the disc and at this point we have $f_{i}=0$ which gives
\begin{equation}   L_{i}= F_{i} - 2 \alpha g^{ab} u_{ai} u_{b}. \end{equation}
We also have
$$   f_{ij} = -L_{ij} + F_{ij} - 2\alpha g^{ab} u_{a ij} u_{b} - 2\alpha
g^{ab} u_{ai} u_{bj}, $$
and at the minumum point $u^{ij} f_{ij} \geq 0$. Hence, at the minimum point,
$$   2 \alpha g^{ab} u_{ab} \leq - u^{ij} L_{ij} + u^{ij} F_{ij}- 2\alpha
g^{ab} L_{a} u_{q}, $$
where we have used the identity
$$ L_{a}= u^{ij} u_{aij}. $$
The defining equation $ u^{ij}_{ij}=-A$ leads to the formula
$$   u^{ij} L_{ij} = u^{ij} L_{i} L_{j} +A, $$
(see \cite{kn:D2},Sect. 2.1) so we get
\begin{equation}   2\alpha g^{ab} u_{ab} \leq A^{-}- u^{ij}L_{i} L_{j}+ u^{ij} F_{ij} - 2\alpha
g^{ab} L_{a} u_{b}. \end{equation}
Now we use (12) to write
$$  u^{ij} L_{i} L_{j} = u^{ij} ( F_{i} - 2\alpha g^{pq} u_{pi} u_{q})( F_{j}-
2\alpha g^{rs} u_{rj} u_{s}) , $$
and expand this out to get
$$  u^{ij} L_{i} L_{j}=  u^{ij}F_{i} F_{j} - 4 \alpha u^{ij} F_{i} g^{rs}
u_{rj} u_{s}+ 4 \alpha^{2} g^{pq} g^{rs} u^{ij} u_{pi} u_{q} u_{rj} u_{s}.
$$ 
This simplifies to 
\begin{equation} u^{ij} L_{i} L_{j}= u^{ij} F_{i} F_{j}- 4\alpha F_{r} u_{s} g^{rs}+ 4 \alpha^{2} g^{pq} g^{rs} u_{q} u_{s} u_{rp}. \end{equation}
Next we use (12) again to write
$$ g^{pq} L_{p} u_{q} = g^{pq}(F_{p} - 2\alpha g^{rs} u_{rp} u_{s}) u_{q},
$$
so
\begin{equation} 4\alpha^{2} g^{pq} g^{rs} u_{q} u_{s} u_{rp} =  2\alpha g^{pq} F_{p} u_{q}
- 2\alpha g^{pq} L_{p} u_{q}. \end{equation}
Combining (13), (14) and (15) we obtain
$$  2 \alpha g^{ab} u_{ab} \leq A^{-} + u^{ij} (F_{ij}-  F_{i} F_{j})
+ 2\alpha F_{r} u_{s} g^{rs}. $$

Now take $F$ to be the function
$F(x)= -2\log( 1-\vert x \vert^{2})$. If $E=e^{-F}$ we have
$$  F_{ij} - F_{i} F_{j}=- E^{-1} E_{ij}$$
and $E= (1-\vert x \vert^{2})^{2}$, so the matrix $(E_{ij})$ is bounded.
Using the formula for the inverse of a $2\times 2$ matrix we get
$$    u^{ij} (F_{ij} - F_{i} F_{j})\leq c \frac{g^{ab} u_{ab}} {(\vert 1-\vert
x \vert^{2})^{2} \det (u_{ij})}, $$
for an easily-computable constant $c$. 
  Similarly the derivative $\nabla F $ is bounded by a multiple of $(1-\vert
  x \vert^{2})^{-1}$ so we obtain, at the minimum point of $f$,
  \begin{equation} 2 \alpha g^{ab} u_{ab} \leq A^{+} + \frac{c}{\det(u_{ij}) (1-\vert x\vert^{2})^{2}}
  g^{ab} u_{ab} + \frac{\alpha c}{ 1-\vert x \vert^{2}}, \end{equation}
  using the fact that $\vert \nabla u \vert \leq 1$. 
  
  Now suppose that, at this minimum point, $$ (1-\vert x \vert^{2})^{2} \det (u_{ij}) \geq \alpha^{-1}.
  $$
  Then we can rearrange to obtain
  $$  \alpha g^{ab} u_{ab} \leq A^{+} + \frac{\alpha c}{(1-\vert x\vert^{2})}.$$
  Since $ 4 \det(u_{ij}) \leq (g^{ab}u_{ab})^{2}$ we have
  $$ (1-\vert x\vert^{2})^{2}  \det(u_{ij})\leq \frac{1}{ 4 \alpha^{2}}\left(
  A^{+}(1-\vert x \vert^{2}) + \alpha c\right)^{2} \leq \left(\frac{ A^{+} +\alpha
  c}{ 2\alpha}\right)^{2}. $$
  So we conclude that, in any event, at the minimum point of $f$,
  $$ (1-\vert x \vert^{2})^{2} \det(u_{ij}) \leq C $$ where
  $$C= \max \left( \alpha^{-1}, \left(\frac{A^{+}
  +\alpha c}{ 2\alpha}\right)^{2}\right). $$
  Taking logarithms, at the minumum point of $f$ we have $-L+F\geq -\log
  C$, so $f\geq -\alpha -\log C$ since $g^{ab} u_{a} u_{b} = \vert \nabla
  u \vert^{2} \leq 1$. So at any point of the disc $-L+F\geq -\alpha -\log
  C$ and in particular when $\vert x\vert \leq \tquart$ we have $\det( u_{ij})\leq \left(\frac{16}{15}\right)^{2}C e^{\alpha}$. This gives our first result, taking any fixed value of $\alpha$.

\

The proof of the second item is very similar. Now we restrict attention to
the open subset $U$ of the unit disc on which $\vert \nabla u \vert <1$ and
 consider the function on $U$
$$f= \log \det(u_{ij}) -\alpha \vert x \vert^{2} +F(\nabla u)$$
where $F$ is a function on the unit disc $\vert \xi \vert < 1$ which tends
to infinity on the boundary. The easiest way to present the proof, in analogy
with preceding case, is to take the Legendre transform $\phi$ of $u$, although
it is not necessary to do so. The point is that the quantity $\det u_{ij}$we want to estimate can also be written as the inverse the determinant of the
Hessian of $\phi$.   We calculate with respect to dual coordinates
$\xi^{i}$. (There is a clash of notation here, in that we would often write
these coordinates with lower indices, to fit in with the previous $x^{i}$,
but that would not be convenient for the calculations we want to perform.) Our function becomes
$$ f= -\log\det \phi_{ij} - \alpha g^{ab} \phi_{a} \phi_{b} +F(\xi), $$
thought of as a function on the unit disc, in $\xi$ coordinates. We write
$L=\log \det \phi_{ij}$, although we should keep in mind that this corresponds
under the Legendre transform to the negative of the function we considered
before. The defining equation for $A$ yields
$$   \phi^{ij} L_{ij} = -A. $$

With these preliminaries in place we can proceed with the argument. At the
minimum we have $L_{i}= F_{i} -2\alpha \phi_{ia} \phi_{b} g^{ab}$ just as
before, and $\phi^{ab} f_{ab} \geq 0$. This leads to
\begin{equation}   2\alpha g^{ab} \phi_{ab}\leq A^{+} + \phi^{ab} F_{ab} -2\alpha g^{ab}
F_{a}
\phi_{b} +4\alpha^{2} \phi_{ia} \phi_{b}\phi_{j} g^{ab} g^{ij}, \end{equation}
(at the minimum point). Now since $\vert \nabla \phi \vert^{2} \leq 1$ we
have
$$  \phi_{ia} \phi_{b} \phi_{j} g^{ab} g^{ij} \leq  \phi_{ab} g^{ab}.$$ 
(To see this, observe that, after rotating coordinates, we can suppose that
$\phi_{2}=0$ at the point in question: then the left hand side of the expression above is $\phi_{11}$
and the right hand side is $\phi_{11}+\phi_{22}$.)
So this time we choose $\alpha<1/4$, in order that the last term in (17) is bounded by $\alpha
g^{ab} \phi_{ab}$, and we obtain
$$ \alpha g^{ab} \phi_{ab} \leq A^{+} +\phi^{ab} F_{ab} - 2\alpha g^{ab}
F_{a} \phi_{b}. $$
We use the same function $F$ as before: $F(\xi)= -2\log(1-\vert \xi\vert^{2})$.
The matrix $F_{ab}$ is bounded by a multiple of $(1-\vert \xi\vert^{2})^{-2}$;
the first derivative $\nabla F$ by a multiple of $(1-\vert \xi\vert^{2})^{-1}$
and the argument proceeds exactly as before.

\end{document}